\documentclass[smallextended,envcountsect]{svjour3}
\usepackage{pifont}
\usepackage{epsfig,color,graphicx,subfigure,amsmath,amssymb,multirow,lscape}
\usepackage[colorlinks,linkcolor=blue,citecolor=blue]{hyperref}
\usepackage{tikz,cite}
\usetikzlibrary{arrows, automata}
\smartqed
\usepackage{graphicx}
\usepackage{amsmath,amssymb}
\journalname{}

\def\pn{\par\smallskip\noindent}

\newtheorem{thm}{Theorem}

\newtheorem{lem}{Lemma}

\newtheorem{rem}{Remark}
\newtheorem{fact}{Fact}
\newtheorem{algo}{Algorithm}
\newtheorem{prob}{Problem}

\def\be{\begin{eqnarray}}
\def\ee{\end{eqnarray}}
\def\ben{\begin{eqnarray*}}
\def\een{\end{eqnarray*}}
\def\ba{\begin{array}}
\def\ea{\end{array}}
\def\bi{\begin{itemize}}
\def\ei{\end{itemize}}
\def\proof {\pn {Proof.} }

\def\endproof{\hfill $\Box$ \vskip .5cm}

\def\cN{{\mathcal N}}
\def\cO{{\mathcal O}}

\def\cS{{\mathcal S}}

\def\bR{{\mathbb R}}
\def\bN{{\mathbb N}}
\def\bB{{\mathbb B}}
\DeclareMathOperator*{\argmin}{argmin}
 \DeclareMathOperator*{\dom}{dom}

\def\prox{{\rm Prox}}

\def\[{\begin{equation}}
\def\]{\end{equation}}

\newcommand{\D}{\Delta}

\renewcommand{\d}{\delta}

\newcommand{\R}{\mathbb R}

\newcommand{\lr}[1]{\left\langle #1\right\rangle}

\begin{document}

\title{First-order primal-dual algorithm with Correction}
\titlerunning{First-order primal-dual algorithm}

\author{Xiaokai Chang$^{1,2}$ \and  Sanyang Liu$^2$   }

\institute{\ding {41} Xiaokai Chang \at xkchang@lut.cn
           \and
1 \ \   School of Science, Lanzhou University of Technology,
               Lanzhou, Gansu, P. R. China. \\
2  \ \  School of Mathematics and Statistics, Xidian University,
               Xi'an, Shaanxi, P. R. China.}

\date{Received: date / Accepted: date}

\maketitle

\begin{abstract}
This paper is devoted to the design of efficient primal-dual algorithm (PDA) for solving convex optimization problems with known saddle-point structure.
We present a new PDA with larger acceptable range of parameters and correction, which result in larger step sizes. The step sizes are predicted by using a local information of the linear operator and corrected by linesearch to satisfy a very weak condition, even weaker than the boundedness of sequence generated.
The convergence and ergodic convergence rate are established for general cases, and in case when one of the prox-functions is strongly convex.
The numerical experiments illustrate the improvements in efficiency from the larger step sizes and  acceptable range of parameters.
\end{abstract}


\section{Introduction}
\label{sec_introduction}
Let $X$, $Y$ be two finite-dimensional real vector spaces equipped with an inner product $\langle\cdot,\cdot\rangle$ and its corresponding norm $\|\cdot\| =\sqrt{\langle\cdot,\cdot\rangle}$. We focus on the following primal problem
\be\label{primal}
\min_{x\in X} f(Kx)+g(x),
\ee
where
\bi
\item $K : X\rightarrow Y$ is a bounded linear operator, with operator norm $L = \|K\|$;
\item $f:\rightarrow(-\infty, +\infty]$ and $g :X \rightarrow(-\infty, +\infty]$ are proper lower semicontinuous convex functions.
\ei

Let $f^*$ denotes the Legendre-Fenchel conjugate of the function $f$, and $K^*$ the adjoint of the operator $K$, then $f^*$ is a proper, convex, lower-semicontinuous (l.s.c.) function. The dual problem of (\ref{primal}) reads as:
\be\label{dual}
\min_{y\in Y} f^*(y)+g^*(-K^*y).
\ee
Actually, problem (\ref{primal}) together with its dual (\ref{dual}) is equivalent to the following convex-concave saddle point problem
\be\label{primal_problem}
\min_{x\in X}\max_{y\in Y} ~~g(x)+\langle Kx,y\rangle-f^*(y).
\ee

By introducing an auxiliary variable $z$, problem (\ref{primal}) can be written as two-block separable convex optimization:
\be\label{two-block}
\min \;\; && f(z)+g(x)\nonumber\\
s.t. \;\;&& Kx-z=0, \\
&& x \in X, ~~z\in Y.\nonumber
\end{eqnarray}
The convex-concave saddle point problem (\ref{primal_problem}) and its primal problem with forms (\ref{primal}) and (\ref{two-block}) are widely presented in many disciplines, including mechanics, signal and image processing, and economics \cite{app1,app2,app3,app4,statistical_learning,11.,FB-Tseng}. Saddle point problems are ubiquitous in optimization as it is a very convenient way to represent many nonsmooth problems, and it in turn often allows to improve the complexity rates from
$\cO(1/\sqrt{N})$ to $\cO(1/N)$. However, the saddle point problem (\ref{primal_problem}) is a typical example where the two simplest iterative methods, the forward-backward method and the Arrow-Hurwicz method \cite{AH}, will not work.

Many efficient methods have been proposed for solving problem (\ref{primal_problem}), for instance, alternating direction method of multipliers (ADMM) \cite{statistical_learning,PC-ADMM,He_ADMM-based,ADMM}, extrapolational gradient methods \cite{6.,extragradient,extragradient-type,13.}, primal-dual algorithms (PDA)  \cite{PC-PDA,CP_PDA,He_PDA,zhang_PDA,inertial-FBF} and their accelerated and generalized versions \cite{Nesterov2004,Acc_PDA,M_PDA}. Here, we concentrate on the most simple first-order PDA and its acceleration. The iterative scheme of the original PDA \cite{CP_PDA,PC-PDA,M_PDA} with fixed step sizes reads as
\be\label{pda_basic}
\left.
\ba{l}
y_{n+1}=\prox_{\sigma f^*}(y_n+\sigma K z_{n}),\\
x_{n+1}=\prox_{\tau g}(x_n-\tau K^*y_{n+1}),\\
z_{n+1}=x_{n+1}+\delta (x_{n+1}-x_n),
\ea\right\}
\ee
where $\delta$ is called an extrapolation parameter, $\tau>0$ and $\sigma> 0$ are regarded as step sizes. When $\delta= 0$ in (\ref{pda_basic}), the primal-dual procedure (\ref{pda_basic}) reduces to the Arrow-Hurwicz
algorithm \cite{AH-PDA}, which has been highlighted in \cite{zhu_TV} for TV image restoration problems. In \cite{CP_PDA}, it was shown that the primal-dual procedure (\ref{pda_basic}) is closely related to many existing methods including the extrapolational gradient method \cite{Popov}, the Douglas-Rachford splitting method \cite{DR-splitting}, and the ADMM.

The convergence of (\ref{pda_basic}) was proved in \cite{CP_PDA} under assumptions $\delta=1$ and $\tau\sigma L^2 < 1$.  Of course, it requires knowing the operator norm $L$ of $K$ to determine step sizes, namely, one need to compute the maximal eigenvalue $\lambda_{\max}(K^*K)$. For simple and sparse $K$ it can not matter, however, for large scale dense matrix $K$ computation of norm is much more expensive. Moreover, the eigenvalues of $K^*K$ can be quite different, so the step sizes governed by the maximal eigenvalue will be very conservative. As first-order algorithm, PDA suffers from slow convergence especially on poorly conditioned problems, when they may take thousands of iterations and still struggle reaching just four digits of accuracy \cite{Acc_PDA}. As a remedy for the slow convergence, diagonal \cite{PC-PDA} and non-diagonal precondition \cite{Acc_PDA} were proposed, numerical efficiency can be improved for some cases, but still there is no strong evidence that such precondition improves or at least does not worsen the speed of convergence of PDA.

Generally, larger step sizes often yield a faster convergence, a linesearch thus was introduced to gain the speed improvement in \cite{M_PDA}. It is known that, the linesearch requires the extra proximal operator or evaluations of $K$, or even both in every linesearch iteration, they will turn to be computationally expensive in situations, where the condition to satisfy is strong and the proximal operator is hard to compute and somewhat expensive. The linesearch in \cite{M_PDA} is to find a proper step size $\tau_n$ satisfying
\be\label{ineq_0}
\sqrt{\beta}\tau_n\|K^*y_{n+1}-K^*y_n\|\leq\alpha\|y_{n+1}-y_n\|
\ee
with $\beta>0$ and $\alpha\in(0,1)$. The important parameter $\alpha$ relating to the step size was restricted on $\alpha \in ]0, 1[$ for guaranteeing the convergence, which will hamper larger step sizes.

\textbf{Contributions.} Our purpose here is to propose an efficient PDA with a prediction-correction procedure (PDA-C for short) to estimate step sizes, rather than using the norm of $K$. The main contribution can be summarized as follows:
\bi
\item We extend the range of $\delta$
to $]\frac{\sqrt{5}-1}{2},+\infty[$, obtain a larger value of $\alpha<\frac{1}{\sqrt{\delta}}$. Namely, $\delta$ can be less than 1 and $\alpha$ can close to $\sqrt{\phi}\approx1.272$  ($\phi=\frac{\sqrt{5}+1}{2}$ is the golden ratio), which will extend the range of step sizes.
\item The step sizes are predicted with low computational cost, by the aid of an inverse local Lipschitz constant of $K^*$, and then corrected when $\delta<1$ to satisfy a very weak condition $\|x_{n+1}-x_n\|<+\infty$. For all the tested problems shown in Section \ref{sec_experiments}, this condition is so weak that the linesearch in Correction step does not start or run only a few times to arrive termination conditions.
\item We prove that PDA-C converges with ergodic rate $\cO(1/N)$ for the primal-dual gap, and introduce an accelerated version of PDA-C under the assumption that the primal or dual problem is strongly convex. The theoretical rate of convergence can be further improved to yield an ergodic rate of convergence $\cO(1/N^2)$ for the primal-dual gap.
\ei

The paper is organized as follows. In Section \ref{sec_preliminarries}, we provide some useful facts and notations. Section \ref{sec_PDAU} is devoted to our basic PDA with correction and updating step sizes. We prove the convergence and establish the ergodic convergence rate for the primal-dual gap. In Section \ref{sec_Acceleration} we propose an accelerated version of PDA-C under the assumption that the primal or dual problem is strongly convex. The implementation and numerical experiments, for solving the LASSO, min-max matrix game and nonnegative least square, are provided in Section \ref{sec_experiments}. We conclude our paper in the final section.

\section{Preliminaries}
\label{sec_preliminarries}
We state the following notations and facts on the well-known properties of the proximal operator and Young's inequality. Some properties are included in textbooks such as \cite{Bauschke2011Convex}.

Let $g :X \rightarrow(-\infty, +\infty]$ be a proper lower semicontinuous convex functions, the proximal operator $\prox_{\lambda g}:X\rightarrow X$  is defined as $\prox_{\lambda g}(x) = (I+\lambda \partial g)^{-1}(x), \lambda>0, x\in X$, and explicitly as
\ben
\prox_{\lambda g}(x)
= \argmin_{y\in X}\left\{ g(y)+\frac{1}{2\lambda}\|x-y\|^2\right\},\quad\forall\, x\in X, \lambda>0.
\een

\begin{fact}\cite{Bauschke2011Convex}\label{fact_proj}
Let $g :X\rightarrow (-\infty, +\infty]$ be a convex function. Then for any $\lambda>0$ and $x\in X$, $p = \prox_{\lambda g}(x)$ if and
only if
\ben
\langle p-x, y-p\rangle\geq \lambda [g(p)-g(y)],~~ \forall y\in X.
\een
\end{fact}

\begin{fact}\label{fact_ab}
Let $\{a_n\}_{n\in\bN}$, $\{b_n\}_{n\in\bN}$ be two nonnegative real sequences and $\exists N>0$ such that
\ben
a_{n+1} \leq a_n-b_n,~~\forall n>N.
\een
Then $\{a_n\}_{n\in\bN}$ is convergent and $\lim_{n\rightarrow \infty} b_n = 0$.
\end{fact}

\begin{fact}\label{fact_Yang}
(Young's inequality)  For any $a, b\geq0$ and $\varepsilon> 0$, we have
$$
ab\leq \frac{a^2}{2\varepsilon} + \frac{\varepsilon b^2}{2}.
$$
\end{fact}

The following identity (cosine rule) appears in many convergence analyses and we will use it many times. For any $x, y, z\in \bR^n$,
\be\label{id}
\langle x-y, x-z\rangle= \frac{1}{2}\|x-y\|^2 +\frac{1}{2}\|x-z\|^2- \frac{1}{2}\|y-z\|^2.
\ee

We assume that the solution set of problem (\ref{primal_problem}) is nonempty and denoted by $\cS$. Let $(\bar{x}, \bar{y})$ be a saddle point of problem (\ref{primal_problem}), i.e. $(\bar{x}, \bar{y})\in\cS$, it therefore satisfies
\ben
K\bar{x} \in \partial f^*(\bar{y}), ~~-(K^*\bar{y})\in \partial g(\bar{x}),
\een
where $\partial f^*$ and $\partial g$ are the subdifferential of the convex functions $f^*$ and $g$. For more details on the theory of saddle point, see \cite{Bauschke2011Convex}. Throughout the paper we will assume that $f$ and $g$ are proper (or simple), in the sense that their resolvent operator has a closed-form representation.

By the definition of saddle point, for any $(\bar{x}, \bar{y})\in\cS$ we have
\be
P_{\bar{x}, \bar{y}}(x) := g(x)-g(\bar{x}) + \langle K^*
\bar{y}, x-\bar{x}\rangle \geq0,~\forall x\in X, \label{P}\\
D_{\bar{x}, \bar{y}}(y) := f^*(y)-f^*(\bar{y}) - \langle K
\bar{x}, y-\bar{y}\rangle \geq 0,~\forall y\in Y. \label{D}
\ee
The primal-dual gap can be expressed as $G_{\bar{x}, \bar{y}}(x,y)= P_{\bar{x}, \bar{y}}(x)+D_{\bar{x}, \bar{y}}(y)$. In certain cases when it is clear which saddle point is considered, we will omit the subscript in $P$, $D$ and $G$. It is also important to highlight that functions $P(\cdot)$ and $D(\cdot)$ are convex for fixed $(\bar{x}, \bar{y})\in\cS$.

\section{Primal-Dual Algorithm with Correction}
\label{sec_PDAU}
In this section, we state our primal-dual algorithm with correction and explore its convergence. The step sizes are predicted by the aid of an inverse local Lipschitz constant of $K^*$ and corrected using linesearch.

\vskip5mm
\hrule\vskip2mm
\begin{algo}
[PDA-C for solving (\ref{primal_problem})]\label{algo1}
{~}\vskip 1pt {\rm
\begin{description}
\item[{\em Step 0.}] Take $\delta\in]\frac{\sqrt{5}-1}{2},+\infty[ $, $\varrho\in]0,1[$ and $1<\nu\leq\mu$. Choose $x_0\in X,$ $y_0\in Y$,  $\lambda_0=\lambda_1>0$, $\beta>0$ and $\alpha \in ]0, \frac{1}{\sqrt{\delta}}[ $. Set $n=0$ and $$\zeta_0=\max\{\|x_0-\prox_{\lambda_0 g}(x_0-\lambda_0 K^*y_0)\|,~~\|y_{0}-\prox_{\lambda_{0} f^*}(y_0+\beta\lambda_{0} K x_{0})\|\}.$$
\item[{\em Step 1.}]
1.a. Compute
\be
 x_{n+1}&=&\prox_{\lambda_n g}(x_n-\lambda_n K^*y_n),\label{x_updating}
\ee
1.b. \textbf{Correct when $\delta<1$}: compute $\zeta_{n+1}=\|x_{n+1}-x_n\|$, check
\be\label{zeta}
\zeta_{n+1}\leq \min\{\mu\zeta_0, ~~\nu\zeta_{n}\},
\ee
if not hold, set $\lambda_n\leftarrow \varrho \lambda_n$, $\lambda_{n+1}\leftarrow \min\{\lambda_n, \lambda_{n+1}\}$ and return to Step 1.a.
\item[{\em Step 2.}] Compute
\be
 z_{n+1}&=&x_{n+1}+\delta (x_{n+1}-x_n),\nonumber \\
 y_{n+1}&=&\prox_{\lambda_{n+1} f^*}(y_n+\beta\lambda_{n+1} K z_{n+1}),\label{y_updating}
 \ee
and update
\be\label{lambda}
\lambda_{n+2}&=& \left\{\ba{cl}
         \min~\left\{{\frac{\alpha\|y_{n+1}-y_{n}\|} {\sqrt{\beta}\|K^*y_{n+1}-K^*y_{n}\|},~~ \lambda_{n+1} }\right\} , &  \mbox{if}\ \ K^*y_{n+1}-K^*y_{n}\neq0, \\
          \lambda_{n+1},& \mbox{otherwise};
          \ea
          \right.
\ee
\item[{\em Step 3.}] Set $n\leftarrow n + 1$ and return to step 1.\\
  \end{description}
}
\end{algo}
\vskip1mm\hrule\vskip5mm

\begin{rem}
For brevity of establishing convergence, different step sizes $\lambda_n$ and $\lambda_{n+1}$ are used in (\ref{x_updating}) and (\ref{y_updating}), respectively. So we have to take two step sizes to compute $x_{n+1}$ and $y_{n+1}$, then obtain the next step size $\lambda_{n+2}$ during each iteration. Furthermore, if $\delta=1$, then $\alpha<1$, in this sense the step sizes agree with that introduced in \cite{CP_PDA,M_PDA}.
\end{rem}

\begin{rem}
Note that the primal and dual variables are symmetrical in the problem (\ref{primal_problem}) and PDA-C, we thus choose a variable with simple proximal operator to correct. In practice, there are many functions with linear (or affine) proximal operator, for instance, $\langle a,x\rangle$, $\frac{1}{2}\|x-a\|^2$ and the indicator function $l_C(x)$ with $C=\{y:\langle a,y\rangle=b\}$ or $C=\bB(c, r)$, a closed ball with a center $c$ and a radius $r> 0$. For these functions, the linesearch becomes extremely simple: it does not require any additional matrix-vector multiplications.
\end{rem}

The aim of Correction step is to bound $\{\|x_{n+1}-x_{n}\|\}$ when $\delta<1$, as convergence analysis requires $\|x_{n+1}-x_{n}\|<+\infty$. From (\ref{zeta}), we have $\zeta_n\leq \nu\zeta_0$ for all $n\geq 1$, and $\zeta_{n+1}\leq\mu \zeta_{n}$ for bounding more tightly due to $\|x_{n+1}-x_n\|\rightarrow 0$. The following lemma shows that the correction procedure described in Algorithm \ref{algo1} is well-defined.

\begin{lem}\label{lem_cor}
The correction procedure always terminates. i.e., $\{\lambda_n\}$ is well defined when $\delta\in]\frac{\sqrt{5}-1}{2}, 1[$.
\end{lem}
\proof
Denote
\ben
A:=\partial g ~~\mbox{and}~~x_{n+1}(\lambda) := \mbox{prox}_{\lambda g}(x_n-\lambda K^*y_n).
\een
From \cite[Theorem 23.47]{Bauschke2011Convex}, we have that $\mbox{prox}_{\lambda g}[x_{n+1}(0)]\rightarrow P_{\overline{\dom A}}[x_{n+1}(0)]$ as $\lambda\rightarrow 0$ ($\overline{\dom A}$ denotes the
closures of $\dom A$), which together with the nonexpansivity of $\mbox{prox}_{\lambda g}$ yields
\ben
&&\|x_{n+1}(\lambda)-P_{\overline{\dom A}}[x_{n+1}(0)]\|\\
&\leq& \|x_{n+1}(\lambda)-\mbox{prox}_{\lambda g}[x_{n+1}(0)]\|+\|\mbox{prox}_{\lambda g}[x_{n+1}(0)] -P_{\overline{\dom A}}[x_{n+1}(0)]\|\\
&\leq&\lambda\|F(y_n)\|+\|\mbox{prox}_{\lambda g}[x_{n+1}(0)] -P_{\overline{\dom A}}[x_{n+1}(0)]\|.
\een
By taking the limit as $\lambda\rightarrow 0$, we deduce that $x_{n+1}(\lambda)\rightarrow P_{\overline{\dom A}}[x_{n+1}(0)]$. Notice that $x_{n+1}(0)=x_n$, we observe $P_{\overline{\dom A}}[x_{n+1}(0)]=x_n$.

By a contradiction, suppose that the correction procedure in Algorithm \ref{algo1} fails to terminate at the $n$-th iteration.  Then, for all $\lambda= \varrho^i \lambda_n$ with $i=0,1,\cdots $, we have $\|x_{n+1}(\lambda)-x_n\| >\min\{\mu\zeta_0, \nu\zeta_{n}\}$.
Since $\varrho^i\rightarrow 0$ as $i\rightarrow\infty$, so $\lambda\rightarrow 0$, this gives a contradiction $0\geq\min\{\mu\zeta_0, \nu\zeta_{n}\}$, which completes the proof.
\endproof

\begin{lem}\label{lem_bound}
Let $\{\lambda_n\}_{n\in\bN}$ be a sequence generated by PDA-C, then $\{\lambda_n\}_{n\in\bN}$ is bounded, $\lim\limits_{n\rightarrow\infty}\lambda_n>0$ and $\lim\limits_{n\rightarrow\infty}\frac{\lambda_n}{\lambda_{n-1}}=1$.
\end{lem}
\proof
First, $\{\lambda_n\}_{n\in\bN}$ is upper bounded. Note that $K^*$ is a $L$-Lipschitz continuous mapping with $L = \|K\|$, we have
\ben
\frac{\alpha\|y_{n+1}-y_{n}\|}{\|K^*y_{n+1}-K^*y_{n}\|}\geq\frac{\alpha\|y_{n+1}-y_{n} \|}{L\|y_{n+1}-y_{n}\|}=\frac{\alpha}{L}
\een
for $K^*y_{n+1}-K^*y_{n}\neq0$. This implies the predicted step $ \{\lambda_n\}_{n\in\bN} $ has a lower bound $\tau=\min\{\frac{\alpha}{\sqrt{\beta}L},\lambda_0\}$, then $\lambda_n\geq\tau$ when $\delta\geq1$. If $\delta<1$, $\{\lambda_n\}$ is well defined from Lemma \ref{lem_cor}, and has a lower bound $\tau=\min\{{\frac{\varrho^{i_0}\alpha}{\sqrt{\beta}L },\lambda_0 }\}$ for some $i_0\geq0$. Notice that  the sequence $\{\lambda_n\}$ is monotonically decreasing, we have
 $\lambda=\lim\limits_{n\rightarrow\infty}\lambda_n>0$, consequently, $\lim\limits_{n\rightarrow\infty}\frac{\lambda_n}{\lambda_{n-1}}=1$.
\endproof

The properties of the generated step sizes, shown in Lemma \ref{lem_bound}, are vital for establishing convergence of PDA-C. In sequel, we give an observation in detail on the convergence by using this properties.

\subsection{Convergence Analysis}
\label{sec_Convergence}
This section devotes to the convergence theorem of PDA-C. First, we give the following lemmas for any $\delta\in]\frac{\sqrt{5}-1}{2},+\infty[$, which play a crucial role in the proof of the main theorem.

\begin{lem}\label{lem1}
Let $(\bar{x}, \bar{y})\in \cS$, and $\{(x_n,y_n)\}_{n\in\bN}$ be a sequence generated by PDA-C. Define
\be\label{eta}
\eta_{n}:=(1+\delta)P(x_n)-\delta P(x_{n-1})+D(y_n),
\ee
then for any  $(x,y)\in X\times Y$, we have
\ben
\lambda_n\eta_{n}
&\leq&\langle x_{n+1}-x_{n}, \bar{x}-x_{n+1}\rangle+\left\langle \frac{1}{\beta}(y_{n}-y_{n-1}), \bar{y}-y_{n}\right\rangle+\left\langle \frac{\lambda_n(z_{n}-x_{n})}{\delta\lambda_{n-1}}, x_{n+1}- z_n\right\rangle \\
&&+ \lambda_n\langle K^*y_{n}-K^*y_{n-1}, z_n-x_{n+1}\rangle.
\een
\end{lem}
\proof
By (\ref{x_updating}), (\ref{y_updating}) and Fact \ref{fact_proj}, we have
\be
\langle x_{n+1}-x_{n} + \lambda_{n}K^*y_{n}, ~~ \bar{x}-x_{n+1}\rangle\geq \lambda_{n}[g(x_{n+1})-g(\bar{x})],\label{temp_01}\\
\left\langle \frac{1}{\beta}(y_{n}-y_{n-1}) - \lambda_{n}K z_{n}, ~~\bar{y}-y_{n}\right\rangle\geq \lambda_{n}[f^*(y_{n})-f^*(\bar{y})].\label{temp_001}
\ee
Similarly as in (\ref{temp_01}), for any $x\in X$ we have
\ben
\langle x_{n}-x_{n-1} + \lambda_{n-1}K^*y_{n-1}, ~~ x-x_{n}\rangle\geq \lambda_{n-1}[g(x_{n})-g(x)].
\een
Substituting $x= x_{n+1}$ and $x= x_{n-1}$ in the inequality above, we obtain
\be
\langle x_{n}-x_{n-1} + \lambda_{n-1}K^*y_{n-1},~~ x_{n+1}-x_{n}\rangle\geq \lambda_{n-1}[g(x_{n})-g(x_{n+1})],\label{tem1}\\
\langle x_{n}-x_{n-1} + \lambda_{n-1}K^*y_{n-1}, ~~ x_{n-1}-x_{n}\rangle\geq \lambda_{n-1}[g(x_{n})-g(x_{n-1})].\label{tem2}
\ee
Multiplying (\ref{tem2}) by $\delta$ and then adding it to (\ref{tem1}) yields
\be\label{tem3}
\langle x_{n}-x_{n-1} + \lambda_{n-1}K^*y_{n-1},~~ x_{n+1}-z_{n}\rangle\geq \lambda_{n-1}[(1+\delta)g(x_{n})-g(x_{n+1})-\delta g(x_{n-1})],
\ee
where we use $z_{n}=x_{n}+\delta (x_{n}-x_{n-1})$. Multiplying (\ref{tem3}) by $\frac{\lambda_{n}}{\lambda_{n-1}}$ and using $z_{n}=x_{n}+\delta (x_{n}-x_{n-1})$ again, we get
\be\label{temp_02}
\left\langle \frac{\lambda_n(z_{n}-x_{n})}{\delta\lambda_{n-1}} + \lambda_{n}K^*y_{n-1}, ~~x_{n+1}-z_{n}\right\rangle\geq \lambda_{n}[(1+\delta)g(x_{n})-g(x_{n+1})-\delta g(x_{n-1})].
\ee
Note that
\be\label{temp_03}
\langle K^*y_{n}-K^*\bar{y}, ~~z_n-\bar{x}\rangle=\langle Kz_n-K\bar{x},~~y_{n}-\bar{y}\rangle,
\ee
adding (\ref{temp_01}) and (\ref{temp_001}) to (\ref{temp_02}) gives
\be
&&\langle x_{n+1}-x_{n}, \bar{x}-x_{n+1}\rangle+\left\langle \frac{1}{\beta}(y_{n}-y_{n-1}), \bar{y}-y_{n}\right\rangle+\left\langle \frac{\lambda_n(z_{n}-x_{n})}{\delta\lambda_{n-1}}, x_{n+1}- z_n\right\rangle \nonumber\\
&&+ \lambda_n\langle K^*y_{n}-K^*y_{n-1},~~ z_n-x_{n+1}\rangle-\lambda_n \langle K^*\bar{y}, ~~ z_n-\bar{x}\rangle +\lambda_{n} \langle K\bar{x}, ~~ y_{n}-\bar{y}\rangle\nonumber\\
&\geq&\lambda_{n}[f^*(y_{n})-f^*(\bar{y})]+\lambda_{n}[(1+\delta)g(x_{n})-g(\bar{x})-\delta g(x_{n-1})].\label{temp_04}
\ee
Recalling the definitions of $P$ and $D$ in (\ref{P}) and (\ref{D}), and
\ben
&&(1+\delta)g(x_{n})-g(\bar{x})-\delta g(x_{n-1})+\langle K^*\bar{y}, z_n-\bar{x}\rangle\\
&=&(1+\delta)[g(x_n)-g(\bar{x}) + \langle K^*
\bar{y}, x_n-\bar{x}\rangle]-\delta[g(x_{n-1})-g(\bar{x}) + \langle K^*
\bar{y}, x_{n-1}-\bar{x}\rangle]\\
&=&(1+\delta)P(x_n)-\delta P(x_{n-1}),
\een
we can rewrite (\ref{temp_04}) as
\ben
&&\langle x_{n+1}-x_{n}, \bar{x}-x_{n+1}\rangle+\left\langle \frac{1}{\beta}(y_{n}-y_{n-1}), \bar{y}-y_{n}\right\rangle+\left\langle \frac{\lambda_n(z_{n}-x_{n})}{\delta\lambda_{n-1}}, x_{n+1}- z_n\right\rangle \\
&&+ \lambda_n\langle K^*y_{n}-K^*y_{n-1}, z_n-x_{n+1}\rangle\\
&\geq&\lambda_n[D(y_n)+(1+\delta)P(x_n)-\delta P(x_{n-1})].
\een
The proof can be completed with the definition of $\eta_{n}$ in (\ref{eta}).
\endproof

\begin{lem}\label{lem2}
Let $(\bar{x}, \bar{y})\in\cS$, and $\{(x_n,y_n)\}_{n\in\bN}$ be a sequence generated by PDA-C. Then, for any $\varepsilon>0$, define
\be
 a_n:&=&\|x_n-\bar{x}\|^2+\frac{1}{\beta}\|y_{n-1}-\bar{y}\|^2+ 2 \lambda_{n-1}(1+\delta)P(x_{n-1}),\label{a_n}\\
 b_n:&=&(\frac{\lambda_{n}}{\delta \lambda_{n-1}}-\frac{\alpha\varepsilon\lambda_{n}}{\lambda_{n+1}})\|x_{n+1}-z_n\|^2+
\left(1-\frac{\lambda_n}{\delta \lambda_{n-1}} \right)\|x_{n+1}-x_n\|^2+\frac{\delta\lambda_n}{ \lambda_{n-1}} \|x_{n}-x_{n-1}\|^2\nonumber\\
&&+\frac{1}{\beta}\left(1- \frac{\alpha\lambda_{n}}{\varepsilon\lambda_{n+1}}\right)\|y_{n}-y_{n-1}\|^2, \label{b_n}
\ee
we have
\ben
a_{n+1}\leq a_{n}-b_n.
\een
\end{lem}
\proof By Lemma \ref{lem1}, (\ref{id}) and Cauchy-Schwarz inequality, we obtain
\be\label{ineq1}
&&\|x_{n+1}-\bar{x}\|^2+\frac{1}{\beta}\|y_{n}-\bar{y}\|^2+ 2 \lambda_n(1+\delta)P(x_n)\nonumber\\
&\leq& \|x_n-\bar{x}\|^2+\frac{1}{\beta}\|y_{n-1}-\bar{y}\|^2 +2\lambda_{n}\delta P(x_{n-1})-2\lambda_n D(y_n)\nonumber\\
&&-\frac{\lambda_{n}}{\delta \lambda_{n-1}}[\|z_n-x_n\|^2+\|x_{n+1}-z_n\|^2]+
\left(\frac{\lambda_n}{\delta \lambda_{n-1}}-1 \right)\|x_{n+1}-x_n\|^2-\frac{1}{\beta}\|y_{n}-y_{n-1}\|^2   \nonumber\\
&&+ 2\lambda_n\|K^*y_{n}-K^*y_{n-1}\|\|x_{n+1}-z_n\|\nonumber\\
&\leq&\|x_n-\bar{x}\|^2+\frac{1}{\beta}\|y_{n-1}-\bar{y}\|^2 +2\lambda_{n}\delta P(x_{n-1})-2\lambda_n D(y_n)\nonumber\\
&&-\frac{\lambda_{n}}{\delta \lambda_{n-1}}[\|z_n-x_n\|^2+\|x_{n+1}-z_n\|^2]+
\left(\frac{\lambda_n}{\delta \lambda_{n-1}}-1 \right)\|x_{n+1}-x_n\|^2  \nonumber\\
&&+ \frac{2\alpha}{\sqrt{\beta}}\frac{\lambda_{n}}{\lambda_{n+1}}\|y_{n}-y_{n-1}\|\|x_{n+1}-z_n\|- \frac{1}{\beta}\|y_{n}-y_{n-1}\|^2.
\ee
Using Fact \ref{fact_Yang} with any $\varepsilon>0$, we have
\be\label{ineq2}
\frac{2}{\sqrt{\beta}}\|y_n-y_{n-1}\|\|z_n-x_{n+1}\|
 &\leq&\frac{1}{\varepsilon \beta}\|y_n-y_{n-1}\|^2+\varepsilon\|x_{n+1}-z_n\|^2.
\ee
Combining (\ref{ineq1}) with (\ref{ineq2}), by $z_n-x_n=\delta(x_{n}-x_{n-1})$ we have
\be\label{ineq5}
&&\|x_{n+1}-\bar{x}\|^2+\frac{1}{\beta}\|y_{n}-\bar{y}\|^2+ 2 \lambda_n(1+\delta)P(x_n)\nonumber\\
&\leq& \|x_n-\bar{x}\|^2+\frac{1}{\beta}\|y_{n-1}-\bar{y}\|^2 +2\lambda_{n}\delta P(x_{n-1})-2\lambda_n D(y_n)\nonumber\\
&&-(\frac{\lambda_{n}}{\delta \lambda_{n-1}}-\frac{\alpha\varepsilon\lambda_{n}}{\lambda_{n+1}})\|x_{n+1}-z_n\|^2-
\left(1-\frac{\lambda_n}{\delta \lambda_{n-1}} \right)\|x_{n+1}-x_n\|^2-\frac{\delta\lambda_n}{ \lambda_{n-1}} \|x_{n}-x_{n-1}\|^2\nonumber\\
&&-\frac{1}{\beta}\left(1- \frac{\alpha\lambda_{n}}{\varepsilon\lambda_{n+1}}\right)\|y_{n}-y_{n-1}\|^2.
\ee
Since $(\bar{x}, \bar{y})$ is a saddle point, then $D(y_{n})\geq0$ and $P(x_{n-1})\geq0$. Together with $\delta>0$ and $0<\delta\lambda_{n}\leq(1+\delta)\lambda_{n-1}$, the proof is completed by the definitions of $a_n$ and $b_n$.
\endproof

Since $\lim\limits_{n\rightarrow\infty}\frac{\lambda_n}{\lambda_{n-1}}=1$, we have $\lim\limits_{n\rightarrow\infty}(1-\frac{\lambda_n}{\delta \lambda_{n-1}})=1-\frac{1}{\delta}\geq0$ for any $\delta\geq1$. But for $\delta\in]\frac{\sqrt{5}-1}{2},1[$, we have $1-\frac{1}{\delta}<0$. So,   convergence of Algorithm \ref{algo1} with $\delta<1$ is different from that with $\delta\geq1$, and hence cannot be established by the similar methods as in \cite{yang,M_PDA,Proximal-extrapolated}. By summing and  integrating terms $\left(1-\frac{\lambda_n}{\delta \lambda_{n-1}} \right)\|x_{n+1}-x_n\|^2$ and $\frac{\delta\lambda_n}{ \lambda_{n-1}} \|x_{n}-x_{n-1}\|^2$, convergence is established when $\delta<1$ using Correction step in the following Theorem.

\begin{thm}\label{thm1}
Let $(\bar{x}, \bar{y})\in\cS$ and $\{(x_n,y_n)\}_{n\in\bN}$ be a sequence generated by PDA-C. Then it is a bounded
sequence in $X\times Y$ and all its cluster points are solutions of (\ref{primal_problem}). Moreover, if $g|\dom g$ is continuous then the whole sequence $\{(x_{n}, y_n)\}$ converges to a solution of (\ref{primal_problem}).
\end{thm}
\proof
By Lemma \ref{lem_bound}, we observe $\lim\limits_{n\rightarrow\infty}\frac{\lambda_n}{\lambda_{n+1}}=1$ and $\lim\limits_{n\rightarrow\infty}\frac{\lambda_n}{\lambda_{n-1}}=1$. Setting $\varepsilon=\frac{1}{\sqrt{\delta}}$ in Lemma \ref{lem2}, since $\alpha<\frac{1}{\sqrt{\delta}}$ and $\delta\in]\frac{\sqrt{5}-1}{2},+\infty[$ we deduce
\be\label{ineq_6}
\left.
\ba{r}
\lim\limits_{n\rightarrow\infty}(\frac{\lambda_{n}}{\delta \lambda_{n-1}}-\frac{\alpha\varepsilon\lambda_{n}}{\lambda_{n+1}})
=\frac{1}{\delta}-\frac{\alpha}{\sqrt{\delta}}>0,\\
\lim\limits_{n\rightarrow\infty}\left(1- \frac{\alpha\lambda_{n}}{\varepsilon\lambda_{n+1}}\right)=1- \alpha \sqrt{\delta}>0,\\
\lim\limits_{n\rightarrow\infty}\left(1-\frac{\delta\lambda_n}{ \lambda_{n-1}} +\frac{\delta\lambda_{n+1}}{ \lambda_{n}}\right)
=1-\frac{1}{\delta}+\delta
>0.
\ea
\right\}
\ee
There exists an integer $N>2,$ such that for any $n>N$,
\be\label{ineq_all}
\left.
\ba{r}
\frac{\lambda_{n}}{\delta \lambda_{n-1}}-\frac{\alpha\varepsilon\lambda_{n}}{\lambda_{n+1}}>0,\\
1- \frac{\alpha\lambda_{n}}{\varepsilon\lambda_{n+1}}>0,\\
1-\frac{1}{\delta}+\frac{\delta\lambda_{n+1}}{ \lambda_{n}}>0,
\ea
\right\}
\ee
which implies $b_n\geq0$ ($n>N$) and $\delta\geq1$. Hence, by $a_n\geq0$, Lemma \ref{lem2} and Fact \ref{fact_ab},  $\{a_n\}_{n\in\bN}$ is convergent and
$\lim\limits_{n\rightarrow\infty} b_n = 0$,  then
\be\label{limit}
\lim\limits_{n\rightarrow\infty}\|x_{n+1}-z_n\|=0, ~~\lim\limits_{n\rightarrow\infty}\|x_{n+1}-x_n\|=0 ~~\mbox{and}~~ \lim\limits_{n\rightarrow\infty}\|y_n-y_{n-1}\|=0.
\ee

Now, let us explore the case $\delta<1$. By the definition of $b_n$ in (\ref{b_n}), for any $M>N+1$, we have
\ben
&&a_{M+1}-a_{N+1}=\sum\limits_{n=N+1}^{M}(a_{n+1}-a_{n})\\
&\leq& -\sum\limits_{n=N+1}^{M}(\frac{\lambda_{n}}{\delta \lambda_{n-1}}-\frac{\alpha\varepsilon\lambda_{n}}{\lambda_{n+1}})\|x_{n+1} -z_n\|^2-\sum\limits_{n=N+1}^{M}\frac{1}{\beta}\left(1- \frac{\alpha\lambda_{n}}{\varepsilon\lambda_{n+1}}\right)\|y_{n}-y_{n-1}\|^2\nonumber\\
&&-\sum\limits_{n=N+2}^{M}\left(1-\frac{\lambda_n}{\delta \lambda_{n-1}} +\frac{\delta\lambda_{n+1}}{ \lambda_{n}} \right)\|x_{n+1}-x_n\|^2\nonumber\\
&&-\frac{\delta\lambda_{N+1}}{ \lambda_{N}} \|x_{N+1}-x_{N}\|^2 +\xi_M\nonumber\\
&\leq&\xi_M,
\een
where $\xi_M=\left(\frac{\lambda_{M+1}}{\delta \lambda_{M}}-1\right)\|x_{M+2}-x_{M+1}\|^2<+\infty $ from Lemma \ref{lem_bound} and Correction step. This together with $a_{n}\geq0$ implies that $\{a_{n}\}_{n\in\bN}$ is  bounded and
\ben
\left.
\ba{r}
0\leq\sum\limits_{n=N+1}^{\infty}(\frac{\lambda_{n}}{\delta \lambda_{n-1}}-\frac{\alpha\varepsilon\lambda_{n}}{\lambda_{n+1}})\|x_{n+1} -z_n\|^2<+\infty,\\
0\leq\sum\limits_{n=N+1}^{\infty}\frac{1}{\beta}\left(1- \frac{\alpha\lambda_{n}}{\varepsilon\lambda_{n+1}}\right)\|y_{n}-y_{n-1}\|^2<+\infty,\\
0\leq\sum\limits_{n=N+2}^{\infty}\left(1-\frac{\lambda_n}{\delta \lambda_{n-1}} +\frac{\delta\lambda_{n+1}}{ \lambda_{n}} \right)\|x_{n+1}-x_n\|^2<+\infty,
\ea\right\}
\een
so (\ref{limit}) is valid as well.

Due to $\|x_n-\bar{x}\|^2\leq a_n$, then $\{x_{n}\}_{n\in\bN}$ is bounded for any $\delta\in]\frac{\sqrt{5}-1}{2},+\infty[$.
By $\|x_{n+1}-x_n\|=\frac{1}{\delta}\|x_{n+1}-y_{n+1}\|$, we obtain
$\lim\limits_{n\rightarrow\infty}\|x_{n}-y_n\|=0$ and then $\{
y_n\}_{n\in\bN}$ is bounded. Let $\{(x_{n_k+1},y_{n_k})\}_{k\in\bN}$ be a subsequence that converges to some cluster $(x^*,y^*)$, then $z_{n_k}\rightarrow x^*$.  Applying  Fact \ref{fact_proj}, we deduce that
\be\label{ineq3}
\left.\ba{l}
\langle x_{n_k+1}-x_{n_k} + \lambda_{n_k}K^*y_{n_k}, ~~ x-x_{n_k+1}\rangle\geq \lambda_{n_k}[g(x_{n_k+1})-g(x)]\\
\left\langle \frac{1}{\beta}(y_{n_k}-y_{n_k-1}) - \lambda_{n_k}K z_{n_k}, ~~y-y_{n_k}\right\rangle\geq \lambda_{n_k}[f^*(y_{n_k})-f^*(y)],
\ea
\right\}
\ee
for any $(x,y)\in X\times Y$, which implies $(x^*,y^*)$ is a saddle problem (\ref{primal_problem}) by passing to the limit and the fact $\lambda_n$ is separated from 0.

We take $(\bar{x},\bar{y})=(x^*,y^*)$ in the definition of $a_n$ and label as $a_n^*$. Notice that $\{\lambda_n\}_{n\in\bN}$ is bounded and $P(x^*, \cdot)$ is continuous when $g|\dom g$ is continuous, hence, $P(x_{n_k})\rightarrow0$ and
\ben
\lim\limits_{n\rightarrow\infty}a_{n+1}^*
&=&\lim\limits_{k\rightarrow\infty}a_{n_{k}+1}^* \\ &=&\lim\limits_{k\rightarrow\infty}\left(\|x_{n_{k}+1}-x^*\|^2+\frac{1} {\beta}\|y_{n_{k}}-y^*\|^2+ 2\lambda_{n_{k}}(1+\delta)P(x_{n_{k}})\right)
= 0,
\een
which means $x_{n+1}\rightarrow x^*$ and $y_{n}\rightarrow y^*$. This completes the proof.
\endproof

\begin{rem}
From \cite{M_PDA}, the condition of $g|\dom g$ to be continuous is not restrictive: it holds when $\dom g$ is an open set (this includes all finite-valued functions) or $g=\delta_C$ for any closed convex set $C$. Moreover, it holds for any separable lower semicontinuous convex function from \cite[Corollary 9.15]{Bauschke2011Convex}.
\end{rem}
\subsection{Ergodic Convergence Rate}
\label{sec_Rate}
In this section, we investigate the convergence rate of the ergodic sequence $\{(X_j, Y_j)\}_{j\in\bN}$ defined in (\ref{definition_XY}). For the case $\delta\geq1$, it can be obtained by the similar technique as that in \cite{M_PDA}, we thus focus on the case when $\delta\in]\frac{\sqrt{5}-1}{2},1[$.

\begin{thm}\label{thm2}
Let $\{(x_n,y_n)\}_{n\in\bN}$ be a sequence generated by PDA-C with $\delta\in]\frac{\sqrt{5}-1}{2},1[$ and $(\bar{x}, \bar{y})\in\cS$. For any $n_1>N$ and $j>n_1$, we define
\be\label{definition_XY}
s_j=\sum_{l=n_1}^j\lambda_l,~~X_j=\frac{\lambda_{n_1}\delta x_{n_1-1} + \sum_{l=n_1}^j\lambda_l z_l}{\lambda_{n_1}\delta + s_j},~~Y_j=\frac{\sum_{l=n_1}^j\lambda_l y_l}{s_j},
\ee
then there exists a sufficient large $J$, when $j>J$ we have
\ben
G(X_j,Y_j) = P(X_j)+D(Y_j)\leq \frac{1}{2s_j}\left[\|x_{n_1}-\bar{x}\|^2+\frac{1}{\beta}\|y_{n_1-1}-\bar{y}\|^2 + 2\delta \lambda_{n_1}P(x_{n_1-1})\right].
\een
\end{thm}
\proof
First of all, combining the definition of $\eta_{n}$ in (\ref{eta}) with the inequality (\ref{ineq5}) yields
\be\label{ineq51}
2\lambda_n\eta_{n}
\leq\|x_n-\bar{x}\|^2-\|x_{n+1}-\bar{x}\|^2+\frac{1}{\beta}\|y_{n-1}-\bar{y}\|^2- \frac{1}{\beta}\|y_{n}-\bar{y}\|^2-b_n,
\ee
where $b_n$ is defined in (\ref{b_n}).
Recalling (\ref{ineq_all}) and summing from $l= n_1$ ($n_1>N$) to $j>n_1$, we get
\ben
&&\|x_{n_1}-\bar{x}\|^2+\frac{1}{\beta}\|y_{n_1-1}-\bar{y}\|^2-\|x_{j+1}-\bar{x}\|^2- \frac{1}{\beta}\|y_{j}-\bar{y}\|^2\\
&&+\left(\frac{\lambda_j}{\delta \lambda_{j-1}}- 1\right)\|x_{j+1}-x_j\|^2-\frac{\delta\lambda_{n_1}}{ \lambda_{n_1-1}} \|x_{n_1}-x_{n_1-1}\|^2\\
&\geq&2\sum_{l=n_1}^j\lambda_l\eta_{l}.
\een
Since $\|x_{n+1}-x_n\|\rightarrow0$ as $ n\rightarrow+\infty$, there exists a sufficiently large $J$ such that for any $j>J$, it holds $\|x_{j+1}-x_j\|\leq\|x_{n_1}-x_{n_1-1}\|\neq0$ (Here we assume $\|x_{n_1}-x_{n_1-1}\|\neq0$). Then
\ben
&&\left(\frac{\lambda_j}{\delta \lambda_{j-1}}- 1\right)\|x_{j+1}-x_j\|^2-\frac{\delta\lambda_{n_1}}{ \lambda_{n_1-1}} \|x_{n_1}-x_{n_1-1}\|^2\\
&\leq&\left(\frac{\lambda_j}{\delta \lambda_{j-1}}- 1-\frac{\delta\lambda_{n_1}}{ \lambda_{n_1-1}}\right)\|x_{j+1}-x_j\|^2\\
&\leq&\left(\frac{1}{\delta}- 1-\frac{\delta\lambda_{n_1}}{ \lambda_{n_1-1}}\right)\|x_{j+1}-x_j\|^2\leq0
\een
from the third inequality of (\ref{ineq_all}) and $\lambda_j\leq\lambda_{j-1}$. Note that
\ben
\sum_{l=n_1}^j\lambda_l\eta_{l}= \lambda_j(1+\delta)P(x_j)+\sum_{l=n_1+1}^j
[(1+\delta)\lambda_{l-1}-\delta\lambda_l)]P(x_{l-1})-\delta \lambda_{n_1}P(x_{n_1-1})+
\sum_{l=n_1}^j\lambda_lD(y_{l}).
\een
By convexity of $P(\cdot)$, we observe
\ben
&&\lambda_j(1+\delta)P(x_j)+\sum_{l=n_1+1}^j
[(1+\delta)\lambda_{l-1}-\delta\lambda_l)]P(x_{l-1})\\
&\geq& (\lambda_{n_1}\delta + s_j)P\left(\frac{\lambda_{n_1}(1 +\delta)x_{n_1} + \sum_{l=n_1+1}^j\lambda_l z_l}{\lambda_{n_1}\delta + s_j}
\right)\\
&=& (\lambda_{n_1}\delta + s_j)P\left(\frac{\lambda_{n_1}\delta x_{n_1-1} + \sum_{l=n_1}^j\lambda_l z_l}{\lambda_{n_1}\delta + s_j}
\right) \geq s_jP(X_j),
\een
where $s_j=\sum_{l=n_1}^j\lambda_l$. Similarly,
\ben
\sum_{l=n_1}^j\lambda_lD(y_l)\geq s_j D\left(\frac{\sum_{l=n_1}^j\lambda_l y_l}{s_j}\right) = s_j D(Y_j).
\een
Hence, we conclude
\ben
\sum_{l=n_1}^j\lambda_l\eta_{l}\geq s_j[P(X_j) + D(Y_j)]- \lambda_{n_1}\delta P(x_{n_1-1}),
\een
and
\ben
G(X_j,Y_j) = P(X_j)+D(Y_j)\leq \frac{1}{2s_j}\left[\|x_{n_1}-\bar{x}\|^2+\frac{1}{\beta}\|y_{n_1-1}-\bar{y}\|^2 + 2\delta \lambda_{n_1}P(x_{n_1-1})\right],
\een
which finishes the proof.
\endproof

Notice that $ \{\lambda_n\}_{n\in\bN}$ has a lower bound $\tau>0$ from the proof of Lemma \ref{lem_bound}. Fix $n_1>N$, we get $s_j\geq (j-n_1+1)\tau$. This implies $s_j\rightarrow +\infty$ when $j\rightarrow+\infty$ and PDA-C has the same $\cO(1/j)$ ergodic rate of convergence when $j>n_1>N$.

\subsection{Heuristics on Nonmonotonic Step Sizes}
\label{sec_improved}
In Algorithm \ref{algo1}, the step size $\{\lambda_n\}_{n\in\bN}$ is updated but in a nonincreasing way, which might be adverse if the algorithm starts in the region with a big curvature of $K^*$. For the purpose of breaking away from overdependence on the few initial step sizes, we choose $\widehat{\lambda}>0$, $n_0>1$ and a sequence $\{\phi_n\}_{n\in\bN}$ with $\phi_n\in[1,\frac{1+\delta} {\delta}]$ and $\phi_n=1$ when $n>n_0$, and update step sizes using following scheme
\be\label{lambda-non}
\lambda_{n+2}&=& \left\{\ba{cl}
         \min~\left\{{\frac{\alpha\|y_{n+1}-y_{n}\|} {\sqrt{\beta}\|K^*y_{n+1}-K^*y_{n}\|},~~ \phi_{n}\lambda_{n+1} },~~\widehat{\lambda}\right\} , &  \mbox{if}\ \ K^*y_{n+1}-K^*y_{n}\neq0, \\
          \lambda_{n+1},& \mbox{otherwise}.
          \ea
          \right.
\ee
Correspondingly, we use $\lambda_{n+1}\leftarrow \min\{\phi_n\lambda_n, \lambda_{n+1}\}$ in Correction step to ensure $\lambda_{n+1}\leq\phi_n\lambda_n$. 

The role of multiplier $\phi_n$ is to allow step sizes to increase, which fulfills that the step sizes can be updated non-monotonically, unlike the updating strategies presented in \cite{yang}. The constant $\widehat{\lambda}$ in Algorithm \ref{algo1} is given only to ensure the upper boundedness of $\{\lambda_n\}$. Hence, it makes sense to choose $\widehat{\lambda}$ quite large.

In this case, the step sizes can be generated non-monotonically when $n<n_0$ but bounded from Lemma \ref{lem_bound}. Consequently, it follows from $\phi_n=1$ when $n\geq n_0$ for given $n_0$ that the sequence $\{\lambda_{n}\}_{n>n_0}$ is monotonically decreasing. This means $\{\lambda_{n}\}$  is convergent,
\ben
\lim_{n\rightarrow\infty}\lambda_{n}>0,~~ \lim_{n\rightarrow\infty}\frac{\lambda_{n}} {\lambda_{n-1}}=1,
\een
and $\left(\frac{\lambda_{n}}{\delta \lambda_{n-1}}-1\right)\|x_{n+1}-x_{n}\|^2<+\infty$. By $\phi_n\in[1,\frac{1+\delta} {\delta}]$ and $\lambda_{n+1}\leq\phi_n\lambda_{n}$, we can deduce $\delta\lambda_{n+1}\leq(1+\delta)\lambda_{n}$. Then Lemmas \ref{lem1} and \ref{lem2} are still valid. Under these conditions, it is not difficult to prove the convergence of Algorithm \ref{algo1} using (\ref{lambda-non}), but its convergence rate is unknown.

\section{Accelerated PDA when $\delta\geq1$.}
\label{sec_Acceleration}
In this section, we consider accelerated version of the primal-dual algorithm when $\delta\geq1$, as the nonnegativity of $\{b_n\}_{n>N}$ can not be ensured for the case $\delta<1$. In many cases, the speed of convergence of PDA crucially depends on the ratio $\beta$ between primal and dual steps. It is shown in \cite{FIST,CP_PDA,M_PDA,Nesterov1983,Nesterov2008} that in case $g$ or $f^*$ are strongly convex, one can modify PDA and derive $\cO(1/N^2)$ convergence by adapting alterable $\beta$. We show that the same holds for PDA-C and when a special strategy is used to update $\lambda_n$. Due to the symmetry of the primal and dual variables in the problem (\ref{primal_problem}) and our method PDA-C, we will only treat the case where $g$ is strongly convex for simplicity, the case where $f^*$ is strongly convex is completely equivalent.

Assume that $g$ is $\gamma$-strongly convex, i.e.,
\ben
g(x_2)-g(x_1)\geq \langle u, x_2-x_1\rangle +
\frac{\gamma}{2}\|x_2-x_1\|^2,~~x_1, x_2\in X, u\in \partial g(x_1),
\een
and the parameter $\gamma$ is known. Exploiting the strong convexity of $g$, we introduce the following accelerated PDA-C (APDA-C).
\vskip5mm
\hrule\vskip2mm
\begin{algo}
[Accelerated PDA-C for solving (\ref{primal_problem}) when $g$ is $\gamma$-strongly convex.]\label{algo2}
{~}\vskip 1pt {\rm
\begin{description}
\item[{\em Step 0.}] Take $\delta\in[1,+\infty[$,  choose $x_0, y_0\in X,$  $\lambda_0=\lambda_1>0$, $\beta_0>0$ and $\alpha \in ]0, \frac{1}{\sqrt{\delta}}[ $. Set $n=0$.
\item[{\em Step 1.}]Compute
\be
 x_{n+1}&=&\prox_{\lambda_n g}(x_n-\lambda_n K^*y_n),\nonumber\\
 z_{n+1}&=&x_{n+1}+\delta (x_{n+1}-x_n),\nonumber \\
 \beta_{n+1}&=&\beta_{n}(1+\gamma \lambda_{n+1}),\label{beta_updating}\\
 y_{n+1}&=&\prox_{\lambda_{n+1} f^*}(y_n+\beta_{n+1}\lambda_{n+1} K z_{n+1}).\nonumber
 \ee
\item[{\em Step 2.}] Update
\be\label{acc_lambda}
\lambda_{n+2}= \left\{\ba{cl}
         \min~\left\{{\frac{\alpha\|y_{n+1}-y_n\|} {\sqrt{\beta_{n+1}}\|K^*y_{n+1}-K^*y_n\|},~~ \frac{\sqrt{\beta_{n}}}{\sqrt{\beta_{n+1}}}\lambda_{n+1} } \right\} , &  \mbox{if}\ \ K^*y_{n+1}-K^*y_n\neq0, \\
          \frac{\sqrt{\beta_{n}}}{\sqrt{\beta_{n+1}}}\lambda_{n+1},& \mbox{otherwise}.
          \ea
          \right.
\ee
\item[{\em Step 3.}] Set $n\leftarrow n + 1$ and return to step 1.\\
  \end{description}
}
\end{algo}
\vskip1mm\hrule\vskip5mm

The main difference of the accelerated variant APDA-C from the basic PDA-C is that now we have to update $\beta_{n+1}$ by $\beta_{n+1}=\beta_{n}(1+\gamma \lambda_{n+1})$ in every iteration, and obtain $\lambda_{n+2}$ from  a special strategy (\ref{acc_lambda}), which will result in the unboundedness  of $\{\beta_{n}\}_{n\in\bN}$ and $\lambda_n\rightarrow0 (n\rightarrow+\infty)$. Even so, the desired properties can be established for the sequences $\{\beta_{n}\}_{n\in\bN}$ and $\{\lambda_{n}\}_{n\in\bN}$, shown in the following Lemma. Also notice that the accelerated algorithm above coincides with PDA-C when a parameter of strong convexity $\gamma= 0$.

\begin{lem}\label{acc_beta}
Let $\{\lambda_n\}_{n\in\bN}$ and $\{\beta_n\}_{n\in\bN}$ be sequences generated by APDA-C, then\\
(i) $0<\lambda_{n}\leq\frac{1+\delta}{\delta}\lambda_{n-1}$, $\lim\limits_{n\rightarrow\infty}\lambda_n=0$ and $\lim\limits_{n\rightarrow\infty} \frac{\lambda_n}{\lambda_{n-1}}=1$;\\
(ii) there exists $C>0$ such that $\beta_n\geq C n^2$ for all $n>0$.
\end{lem}
\proof
(i) The result $0<\lambda_{n}\leq\frac{1+\delta}{\delta}\lambda_{n-1}$ is clear as $\frac{\sqrt{\beta_{n}}}{\sqrt{\beta_{n+1}}}<1$ and $\delta\geq1$. Let $\sigma_n:=\sqrt{\beta_{n-1}} \lambda_n$, using (\ref{acc_lambda}) yields
\ben
\sigma_{n+2}= \left\{\ba{cl}
         \min~\left\{{\frac{\alpha\|y_{n+1}-y_n\|} {\|K^*y_{n+1}-K^*y_n\|},~~ \sigma_{n+1} }\right\} , &  \mbox{if}\ \ K^*y_{n+1}-K^*y_n\neq0, \\
          \sigma_{n+1},& \mbox{otherwise}.
          \ea
          \right.
\een
By the similar techniques as in Lemma \ref{lem_bound}, $\{\sigma_n\}_{n\geq1}$ is bound and has a lower bound $\min\{{\frac{\alpha}{L },\sigma_1 }\}$, then  its limit $\sigma=\lim\limits_{n\rightarrow\infty}\sigma_n$ exists and $\sigma>0$.

Suppose that $\lambda_{n}\nrightarrow0$ when $n\rightarrow+\infty$, we observe $\beta_n\rightarrow+\infty$ as $\beta_{n+1}=\beta_{n}(1+\gamma \lambda_{n+1})$ and $\lambda_{n}>0$, then $\sigma_n\rightarrow+\infty$, which is a contradiction. Thus, we have $\lambda_{n}\rightarrow0$. Consequently, we deduce
\ben
\lim\limits_{n\rightarrow\infty}\frac{\beta_{n+1}}{\beta_{n}}= \lim\limits_{n\rightarrow\infty}(1+\gamma \lambda_{n+1})=1,
\een
and
\ben
\lim\limits_{n\rightarrow\infty}\frac{\lambda_{n+1}}{\lambda_{n}}= \lim\limits_{n\rightarrow\infty}\frac{\sigma_{n+1}}{\sigma_{n}} \frac{\sqrt{\beta_{n-1}}}{\sqrt{\beta_{n}}}=1.
\een
(ii) By (\ref{acc_lambda}) and $\phi_n\geq1$, we find
\ben
\lambda_{n+2}= \frac{\alpha\|y_{n+1}-y_n\|} {\sqrt{\beta_{n+1}}\|K^*y_{n+1}-K^*y_n\|}\geq \frac{\alpha} {\sqrt{\beta_{n+1}}L}
\een
or
\ben
\lambda_{n+2}\geq \frac{\sqrt{\beta_{n}}}{\sqrt{\beta_{n+1}}}\lambda_{n+1}\geq \frac{\sqrt{\beta_{n}}}{\sqrt{\beta_{n+1}}} \frac{\sqrt{\beta_{n-1}}}{\sqrt{\beta_{n}}}\cdots \frac{\sqrt{\beta_{n_0}}}{\sqrt{\beta_{n_0+1}}}\frac{\alpha} {\sqrt{\beta_{n_0}}L}=\frac{\alpha} {\sqrt{\beta_{n+1}}L}
\een
for some $n_0<n$ such that $\lambda_{n_0+2}= \frac{\alpha\|y_{n_0+1}-y_{n_0}\|} {\sqrt{\beta_{n_0+1}}\|K^*y_{n_0+1}-K^*y_{n_0}\|}$. Then, we have
\be\label{beta_ineq}
\beta_{n+1}=\beta_{n}(1+\gamma \lambda_{n+1})\geq \beta_{n}\left(1+\gamma \frac{\alpha}{\sqrt{\beta_{n}}L}\right) =\beta_{n}+\gamma \frac{\alpha\sqrt{\beta_{n}}}{L}.
\ee
By induction, there exists $C>0$ such that $\beta_{n}\geq Cn^2$ for all $n>0$.
\endproof

From Lemma \ref{acc_beta} (i), Lemmas \ref{lem1} and \ref{lem2} are still valid with $\beta_n$ instead of $\beta$, but Theorem \ref{thm1} is not necessarily in place due to $\lambda_n\rightarrow0$. In sequel, we explain that our accelerated method APDA-C yields essentially the same rate $\cO(1/j^2)$ of convergence for the primal dual gap, though the convergence of $\{y_n\}_{n\in\bN}$ is not able to prove.

Instead of (\ref{temp_01}), now one can use a stronger inequality
\be\label{stronger_temp1}
\langle x_{n+1}-x_{n} + \lambda_{n}K^*y_{n}, ~~ \bar{x}-x_{n+1}\rangle\geq \lambda_{n}[g(x_{n+1})-g(\bar{x})+
\frac{\gamma}{2}\|x_{n+1}-\bar{x}\|^2].
\ee
In turn, using (\ref{stronger_temp1}) and the definition of $\eta_n$ yields a stronger version of (\ref{ineq51}) (also with $\beta_n$ instead of $\beta$).
\be\label{stronger_ineq51}
&&(1+\gamma\lambda_{n+1})\|x_{n+1}-\bar{x}\|^2+\frac{\beta_{n+1}}{\beta_{n}} \frac{1}{\beta_{n+1}}\|y_{n}-\bar{y}\|^2+ 2 \lambda_{n}\eta_{n}\nonumber\\
&\leq& \|x_n-\bar{x}\|^2+ \frac{1}{\beta_{n}}\|y_{n-1}-\bar{y}\|^2 -b_n.
\ee

Since $\frac{\beta_{n+1}}{\beta_{n}}=1+\gamma \lambda_{n+1}$ in APDA-C. For brevity let
\ben
A_{n}:=\|x_n-\bar{x}\|^2+ \frac{1}{\beta_{n}}\|y_{n-1}-\bar{y}\|^2,
\een
then (\ref{stronger_ineq51}) gives
\ben
\beta_{n+1}A_{n+1}+2 \beta_{n}\lambda_{n}\eta_{n}
\leq \beta_{n}A_{n}-\beta_{n}b_n.
\een
Recalling (\ref{ineq_all}), for any $n>N$, we have
\ben
\beta_{n+1}A_{n+1}+2 \beta_{n}\lambda_{n}\eta_{n}\leq \beta_{n}A_{n}.
\een
Fix $n>N$ and sum from $l= n$ to $j>n$, we get
\ben
\beta_{n}A_{n}-\beta_{j+1}A_{j+1}
\geq2\sum_{l=n}^j\beta_{l}\lambda_l\eta_{l}.
\een

Defining
\ben
s_j=\sum_{l=n}^j\beta_l\lambda_l,~~X_j=\frac{\lambda_n\delta x_{n-1} + \sum_{l=n}^j\beta_l\lambda_l z_l}{\beta_n\lambda_n\delta + s_j},~~Y_j=\frac{\sum_{l=n}^j\beta_l\lambda_l y_l}{s_j},
\een
for any $j>n$, and using the similar process as in Theorem \ref{thm2}, we observe
\ben
\beta_{j+1}A_{j+1}+2s_j G(X_j,Y_j) \leq \beta_{n}A_n+ 2\delta \beta_{n}\lambda_nP(x_{n-1}).
\een
From this we deduce that the sequence $\{y_n\}_{n\in\bN}$ is bounded and
\ben
G(X_j,Y_j) &\leq& \frac{1}{2s_j} [\beta_{n}A_n+ 2\delta\beta_{n} \lambda_nP(x_{n-1})],\\
\|x_{j+1}-\bar{x}\|^2&\leq&  \frac{1}{2\beta_{j+1}}\left[\beta_{n}A_n+2 \delta \beta_{n}\lambda_nP(x_{n-1})\right].
\een

Then follows Lemma \ref{acc_beta}, for some constant $C_1 > 0$ we have
\ben
\|x_{j+1}-\bar{x}\|^2\leq  \frac{C_1}{(j+1)^2}.
\een
From (\ref{beta_ineq}), we get $\beta_{n+1}-\beta_{n}\geq \gamma \frac{\alpha\sqrt{\beta_{n}}}{L}$. Since
\ben \beta_n\lambda_{n}=\beta_n\lambda_{n+1}\frac{\lambda_{n}}{\lambda_{n+1}} =\frac{\lambda_{n}}{\lambda_{n+1}} \frac{\beta_{n+1}-\beta_{n}}{\gamma}\geq \frac{\lambda_{n}}{\lambda_{n+1}} \frac{\alpha\sqrt{\beta_{n}}}{L}\geq \frac{\lambda_{n}}{\lambda_{n+1}} \frac{\alpha\sqrt{C}n}{L},
\een
we obtain $s_j=\sum_{l=n}^j\beta_l\lambda_l=\cO(j^2)$ using $\lim\limits_{n\rightarrow\infty}\frac{\lambda_{n}}{\lambda_{n+1}}=1$. This means that for some constant $C_2 > 0$,
\ben
G(X_j, Y_j)\leq \frac{C_2}{j^2}.
\een
Finally, we have shown the following result:
\begin{thm}
Let $\{(x_n,y_n)\}_{n\in\bN}$ be a sequence generated by APDA-C. Then $\|x_j-\bar{x}\|=\cO(1/j)$ and $G(X_j, Y_j) = \cO(1/j^2)$.
\end{thm}

\section{Numerical Experiments}
\label{sec_experiments}
We present numerical results to demonstrate the computational performance of  PDA-C (Algorithm \ref{algo1} using (\ref{lambda-non}) to update step sizes) and its acceleration (Algorithm \ref{algo2}) \footnote{All codes are available at http://www.escience.cn/people/changxiaokai/Codes.html} for solving  some minimization problems with saddle-point structure. The following state-of-the-art algorithms are compared to  investigate the computational efficiency:
\begin{itemize}
\item  Tseng's forward-backward-forward splitting method used as in  \cite[Section 4]{Proximal-extrapolated} (denoted by ``FBF"), with $\beta= 0.7, \theta= 0.99$;
\item Proximal gradient method (denoted by ``PGM"), with fixed step $\frac{1}{\|K\|^2}$;
\item Proximal extrapolated gradient methods \cite[Algorithm 2]{Proximal-extrapolated} (denoted by ``PEGM"), with line search and $\alpha= 0.41, \sigma= 0.7$;
\item Primal-Dual algorithm with linesearch \cite{M_PDA} (denoted by ``PDA-L"), with $\mu = 0.7$, $\alpha= 0.99$, and $\tau_0=\frac{\sqrt{\min\{n,m\}}}{\|K\|_F}$.
\item FISTA \cite{FIST,Nesterov1983} with standard linesearch (denoted by ``FISTA"), with $\beta= 0.7, \lambda_0= 1$;
\end{itemize}

We denote the random number generator by $seed$ for generating data again in Python 3.8. All experiments are performed on an Intel(R) Core(TM) i5-4590 CPU@ 3.30 GHz PC with 8GB of RAM running on 64-bit Windows operating system.

There are many choices of the sequence $\{\phi_n\}_{n\in\bN}$, but in the earlier iterations the large range of $\lambda_n$ is benefit for selecting proper step size, we thus use
\be\label{phi0}
\phi_n=\left\{
\ba{rl}\frac{1+\delta}{\delta},~&~\mbox{if}~n\leq\hat{n};\\
\frac{1+\delta+n-\widehat{n}}{\delta+n-\widehat{n}},~&~\mbox{if}~n>\widehat{n};\\
1,~&~\mbox{if}~n>n_0,
\ea
\right.
\ee
for given $\hat{n}, n_0\in\bN$. For PDA-C, we set $\varrho=0.7$, $\alpha=1.27$ and $\delta=0.62$ unless otherwise stated. For APDA-C, we set $\alpha=0.99$ and $\delta=1$.

\begin{prob}[LASSO]\label{pro_1}
We want to minimize:
$$
\min_x \phi(x):=\frac 1 2 ||Kx-b||^2 + \mu  ||x||_1
$$
where $K\in \R^{m\times n}$ is a matrix data, $b\in \R^m$ is a given observation, and $x\in \R^n$ is an unknown signal.
\end{prob}

We can rewrite the problem above in a primal-dual form as follows:
\be\label{pro_11}
\min_{x\in \R^n} \max_{y\in \R^m} g(x)+\langle Kx,y\rangle-f^*(y),
\ee
where $f(p) = \frac 1 2 ||p-b||^2$, $f^*(y) = \frac 1 2 ||y||^2 + (b,y) = \frac 1 2 ||y+b||^2 -\frac{1}{2}||b||^2$ and $g(x) = \mu ||x||_1$.

We set $seed=1$ and generate some random $w\in \bR^n$ in which $s$ random coordinates are drawn from $\cN(0,1)$ and the rest are zeros. Then we generate $\nu\in \bR^m$ with entries drawn from $\cN(0, 0.1)$ and set $b = Kw +\nu$.  The matrix $K\in \bR^{m\times n}$ is constructed in one of the following ways:
\bi
\item[1.] $n = 1000$, $m = 200$, $s = 10$, $\mu= 0.1$. All entries of $K$ are generated independently from $\cN(0; 1)$. The $s$ entries of $w$ are drawn from the uniform distribution in $[-10, 10]$.
\item[2.] $n = 2000$, $m = 1000$, $s = 100$, $\mu= 0.1$. All entries of $K$ are generated independently from $\cN(0, 1)$. The $s$ entries of $w$ are drawn from $\cN(0, 1)$.
\ei

For the primal-dual form (\ref{pro_11}) of Problem \ref{pro_1}, we apply primal-dual methods, and for the problem in a primal form we apply PGM and FISTA. For this we set $h(x) = f(Ax)$ and get $\nabla h(x) = K^*(Kx-b)$. The values of parameters are set as in \cite{M_PDA}, here we rewrite them to facilitate the readers. For PGM and FISTA a fixed step size $\alpha= \frac{1}{\|K\|^2}$ is used. For PDA (\ref{pda_basic}) we use $\sigma= \frac{1}{20\|K\|}, ~\tau=\frac{20}{\|K\|}$. For PDA-L and PDA-C we set $\beta= \frac{1}{400}$ and for PDA-C we set $\hat{n}=5000$. The initial points for all methods are $x_0 = (0,\cdots, 0)$ and $y_0=Kx_0-b$.

\begin{figure}[htp]
\centering
\subfigure[$\phi(x)-\phi^*$]{
\includegraphics[width=0.45\textwidth]{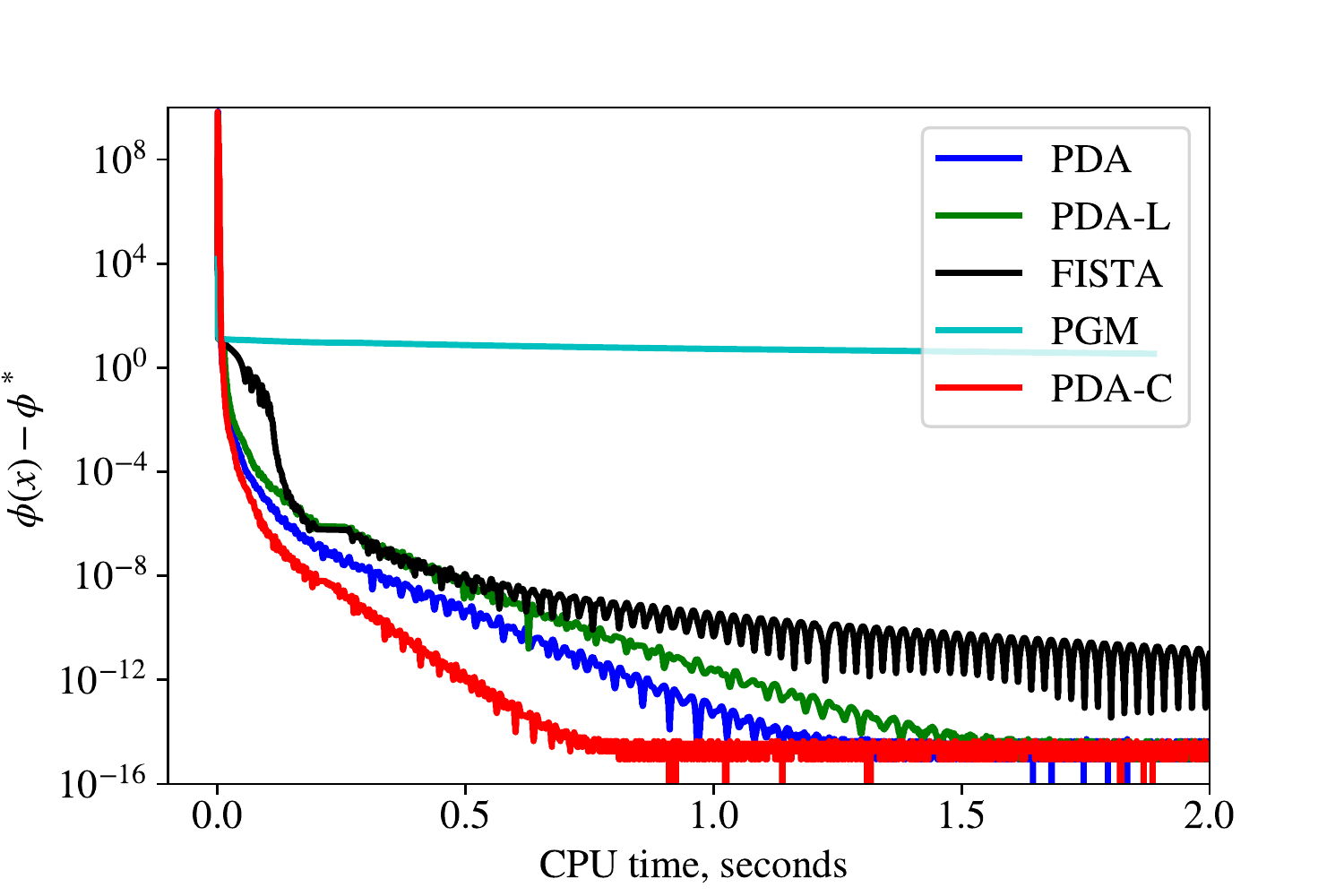}}
\subfigure[$\lambda_n$ (or $\tau_k$)]{
\includegraphics[width=0.45\textwidth]{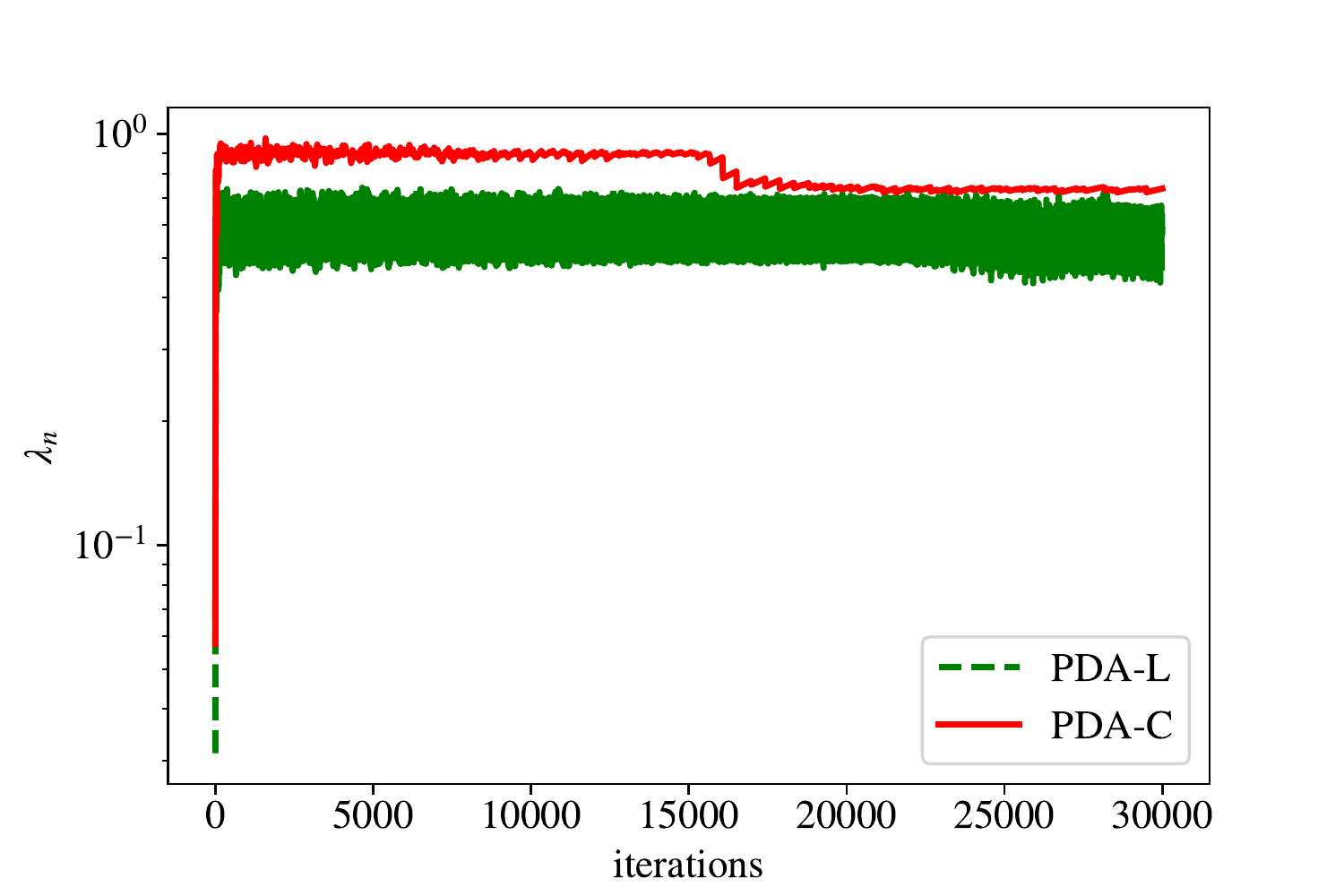}}
\caption{Comparison of $\phi(x)-\phi^*$ and $\lambda_n$ (or $\tau_k$) for solving Problem \ref{pro_1} generated by the first way. }
\label{Fig 1} 
\end{figure}

\begin{figure}[htp]
\centering
\subfigure[$\phi(x)-\phi^*$.]{
\includegraphics[width=0.45\textwidth]{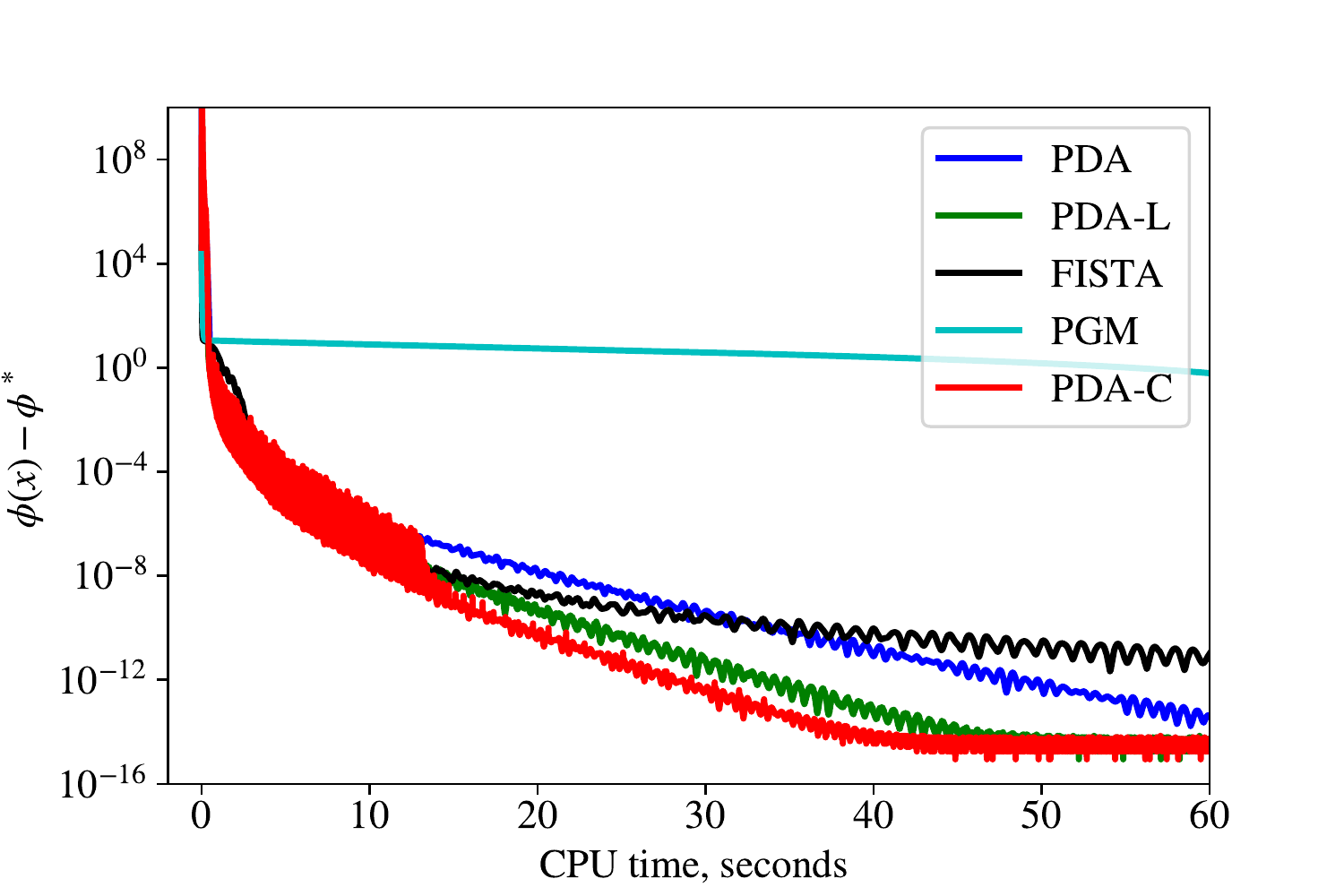}}
\subfigure[$\lambda_n$ (or $\tau_k$)]{
\includegraphics[width=0.45\textwidth]{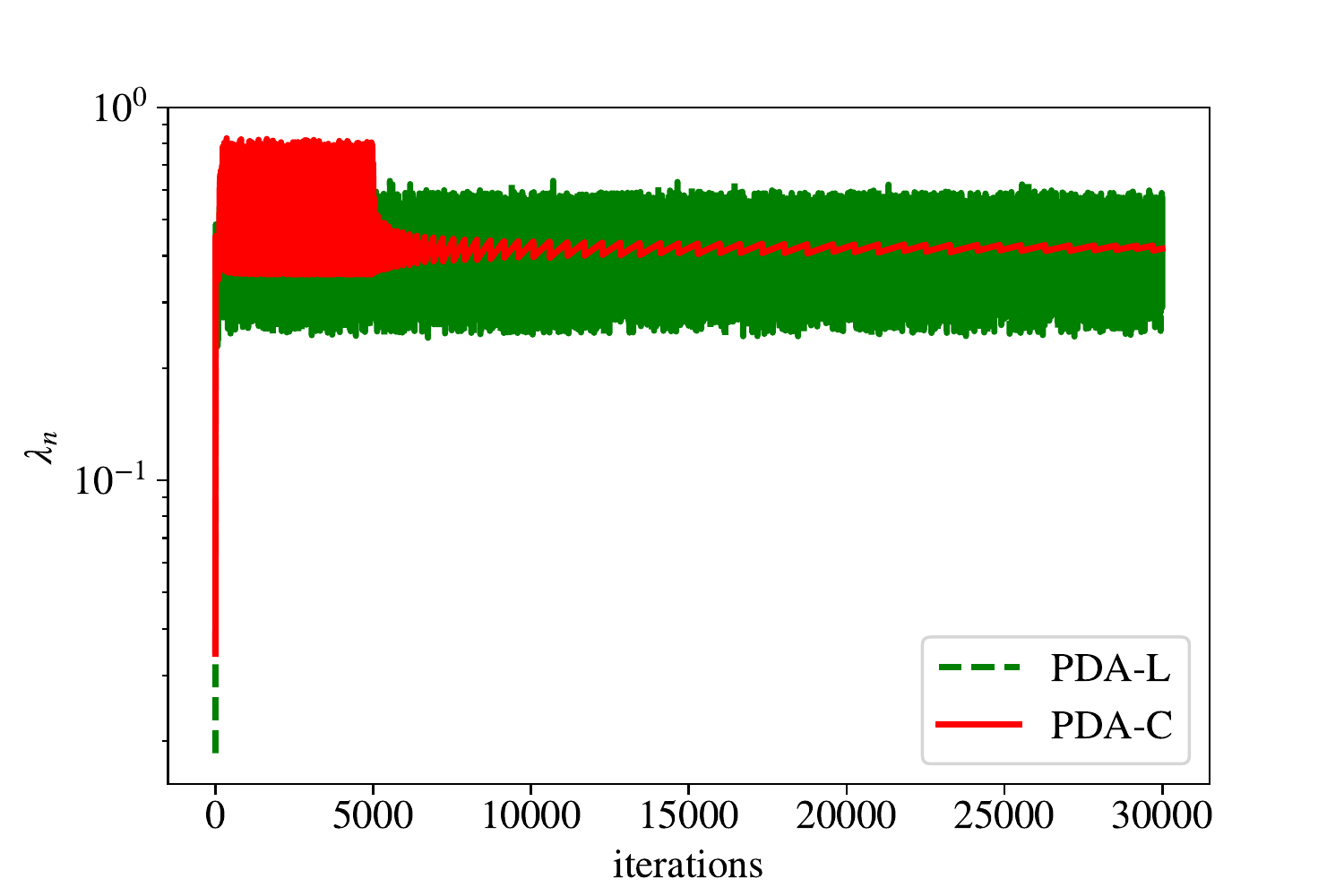}}
\caption{Comparison of $\phi(x)-\phi^*$ and $\lambda_n$ (or $\tau_k$) for solving Problem \ref{pro_1} generated by the second way. }
\label{Fig 2} 
\end{figure}

To illustrate how does the values  $\phi(x_n)-\phi^*$ ($\phi^*=\phi(\bar{x})$ with $(\bar{x},\bar{y})\in\cS$) and $\lambda_n$ for PDA-C (or $\tau_k$ for PDA-L) change over iterations, we give two convergence plots for the maximum number of iterations set at 30,000. From the results shown in Fig. \ref{Fig 1} and \ref{Fig 2}, primal-dual methods show better performance for the instances of Problem \ref{pro_1}, though they require a tuning the parameter $\beta$. 

For the tested problems, PDA-C needs to correct less that 20 times linesearch, so PDA-C with $\delta=0.62$ is more efficient than PDA-L. The advantage of PDA-C is a larger interval
for possible step size $\lambda_n$, see Fig.\ref{Fig 1} (b) and Fig.\ref{Fig 2} (b), which resulted from the smaller choice of $\delta$ and the larger value of $\alpha$.

\begin{prob}[Min-max matrix game]\label{pro_3}
The second problem is the following min-max matrix game
\be\label{mm_pro}
\min_{x \in \D_n}\max_{y\in \D_m} \lr{Kx, y},
\ee
where $x\in \R^n$, $y\in \R^m$, $K\in \R^{m\times n}$, and $\Delta_m$,
$\D_n$ denote the standard unit simplices in $\R^m$ and $\R^n$
respectively.
\end{prob}

The variational inequality formulation of (\ref{mm_pro}) is:
$$\lr{F(z^*),z-z^*} + G(z) - G(z^*) \geq 0 \quad \forall z \in Z,$$
where
$$
Z = \R^n\times \R^m,\quad z=\binom{x}{y},\quad F = \begin{bmatrix} 0 & K^*\\ -K & 0\end{bmatrix}, \quad G(z) = \d_{\D_n}(x) + \d_{\D_m}(y).
$$

For a comparison we use a primal-dual gap (PD gap) as in \cite{M_PDA}, which can be easily computed for a feasible pair $(x, y)$, defined as
\ben
G(x, y) := \max_i(Kx)_i - \min_j(K^*y)_j.
\een
We use an auxiliary point (see \cite{modified-FB}) to compute the primal-dual gap for Tseng's method FBF, as its iterates may be infeasible. The initial point in all cases is chosen as $x_0 =\frac{1}{n}(1,\cdots,1)$ and $y_0 =\frac{1}{m}(1,\cdots,1)$.
We use the algorithm from \cite{projections} to compute projection onto the unit simplex. For PDA we set $\tau= \sigma= \frac{1}{\|K\|}$, which we compute in advance.
The input data for FBF and PEGM with linesearch are taken the same as in \cite{Proximal-extrapolated}. For PDA-L and PDA-C we set $\beta= 1$ (the same as $\tau=\sigma$ in PDA) and for PDA-C we set $\hat{n}=40,000$. Since PDA-C with $\delta=1$ performs well than smaller $\delta$ for the min-max matrix game, we apply PDA-C with $\delta=1$ to testify.

We consider four differently generated samples of the matrix $K\in \bR^{m\times n}$ with $seed=100$ as in \cite{M_PDA}:
\bi
\item[1.] $m = 100$, $n = 100$. All entries of $K$ are generated independently from the uniform distribution in $[-1, 1]$;
\item[2.] $m = 100$, $n = 100$. All entries of $K$ are generated independently from the normal distribution $\cN (0,1)$;
\item[3.] $m = 500$, $n = 100$. All entries of $K$ are generated independently from the normal distribution $\cN (0,10)$;
\item[4.] $m = 100$, $n = 200$. All entries of $K$ are generated independently from the uniform distribution in $[0, 1]$.
\ei

\begin{figure}[htp]
\centering
\subfigure[Example 1.]{
\includegraphics[width=0.45\textwidth]{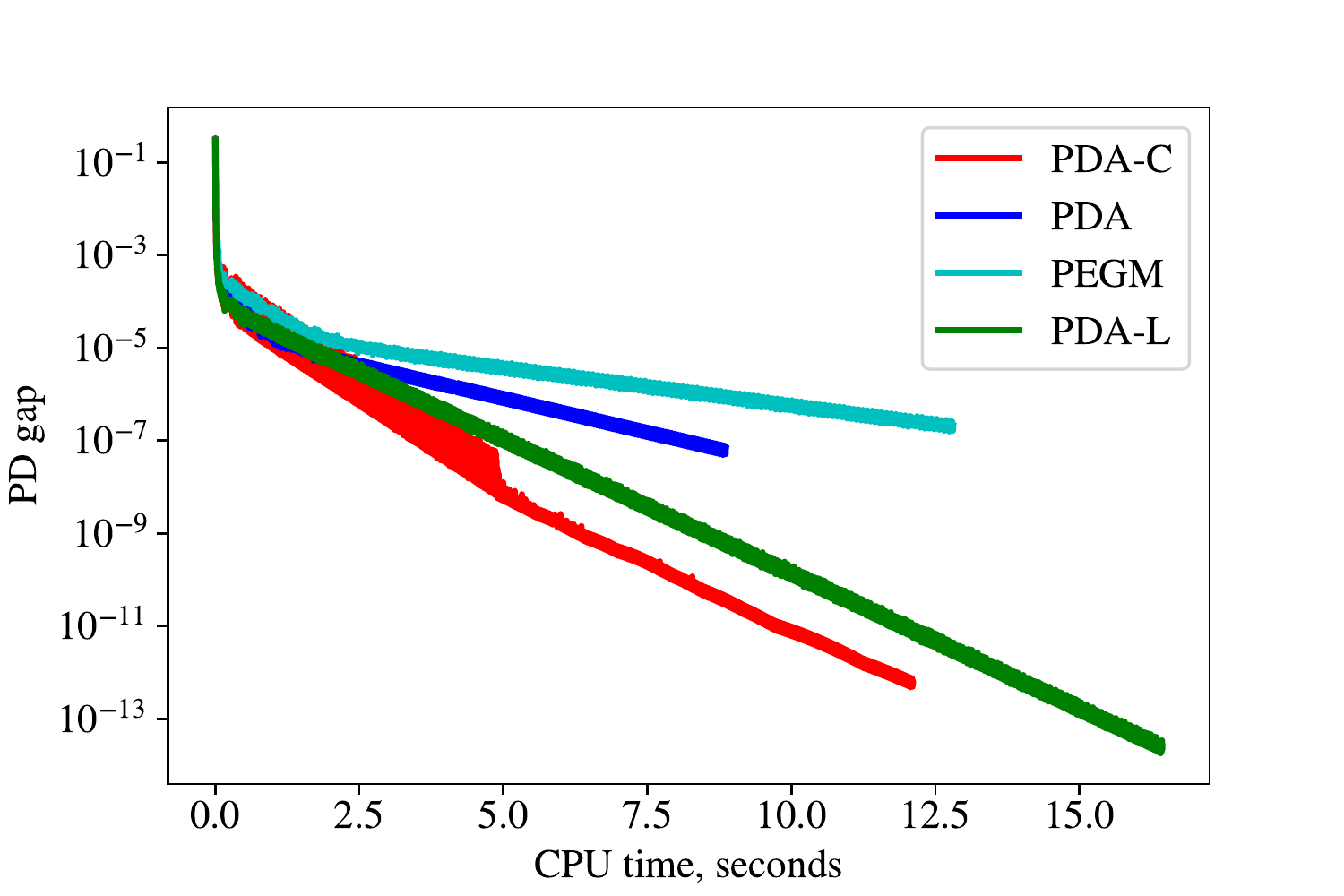}}
\subfigure[Example 2.]{
\includegraphics[width=0.45\textwidth]{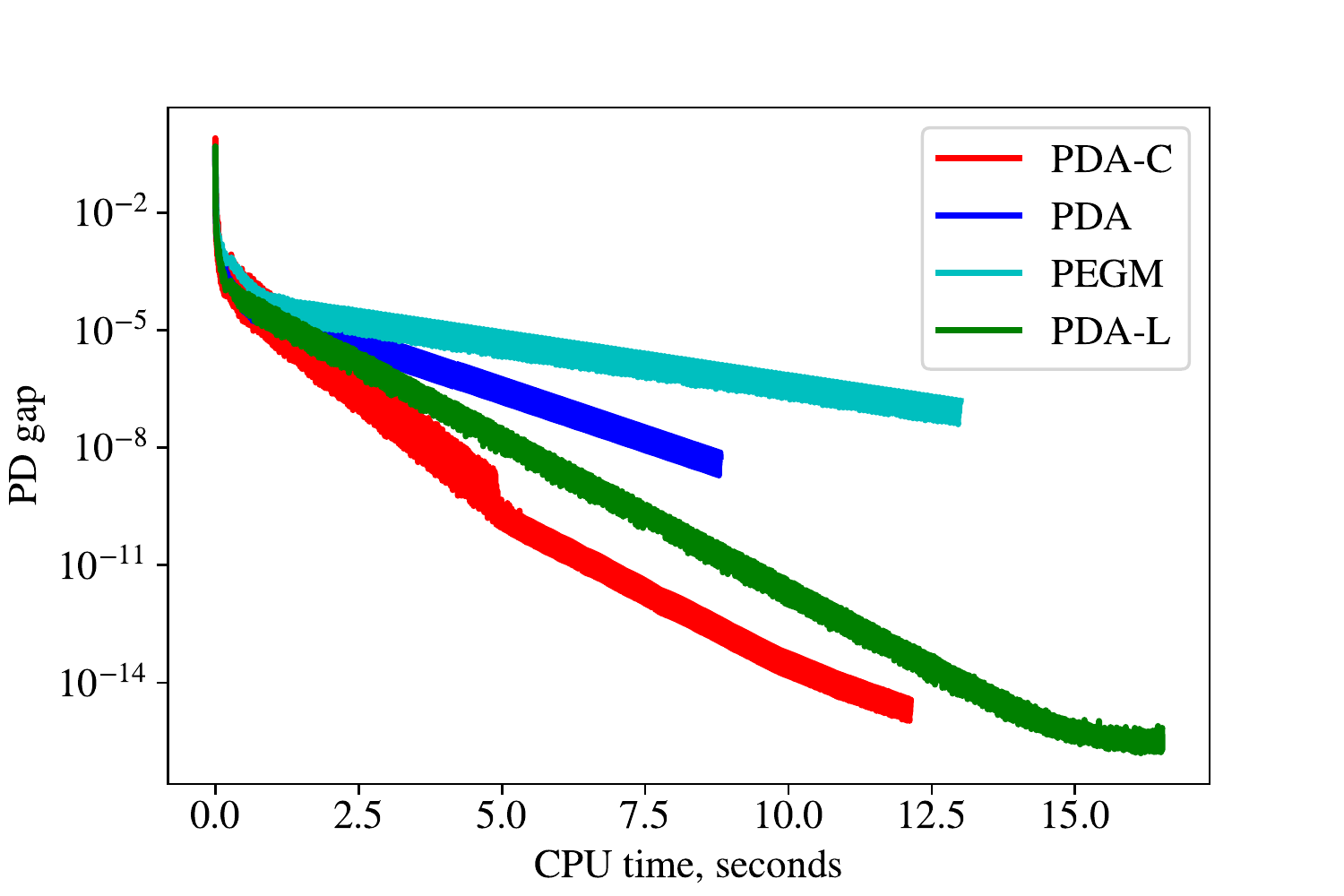}}
\subfigure[Example 3.]{
\includegraphics[width=0.45\textwidth]{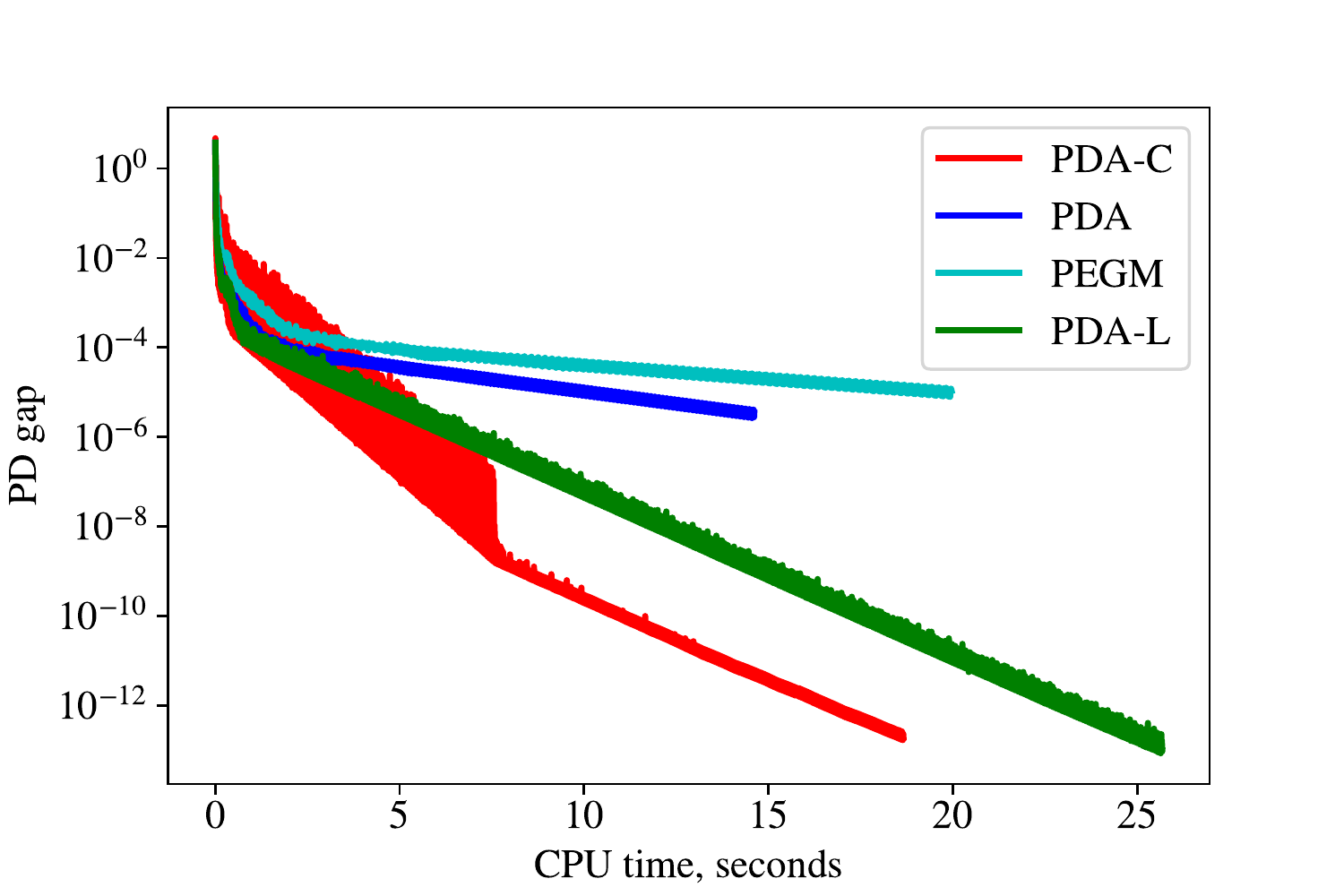}}
\subfigure[Example 4.]{
\includegraphics[width=0.45\textwidth]{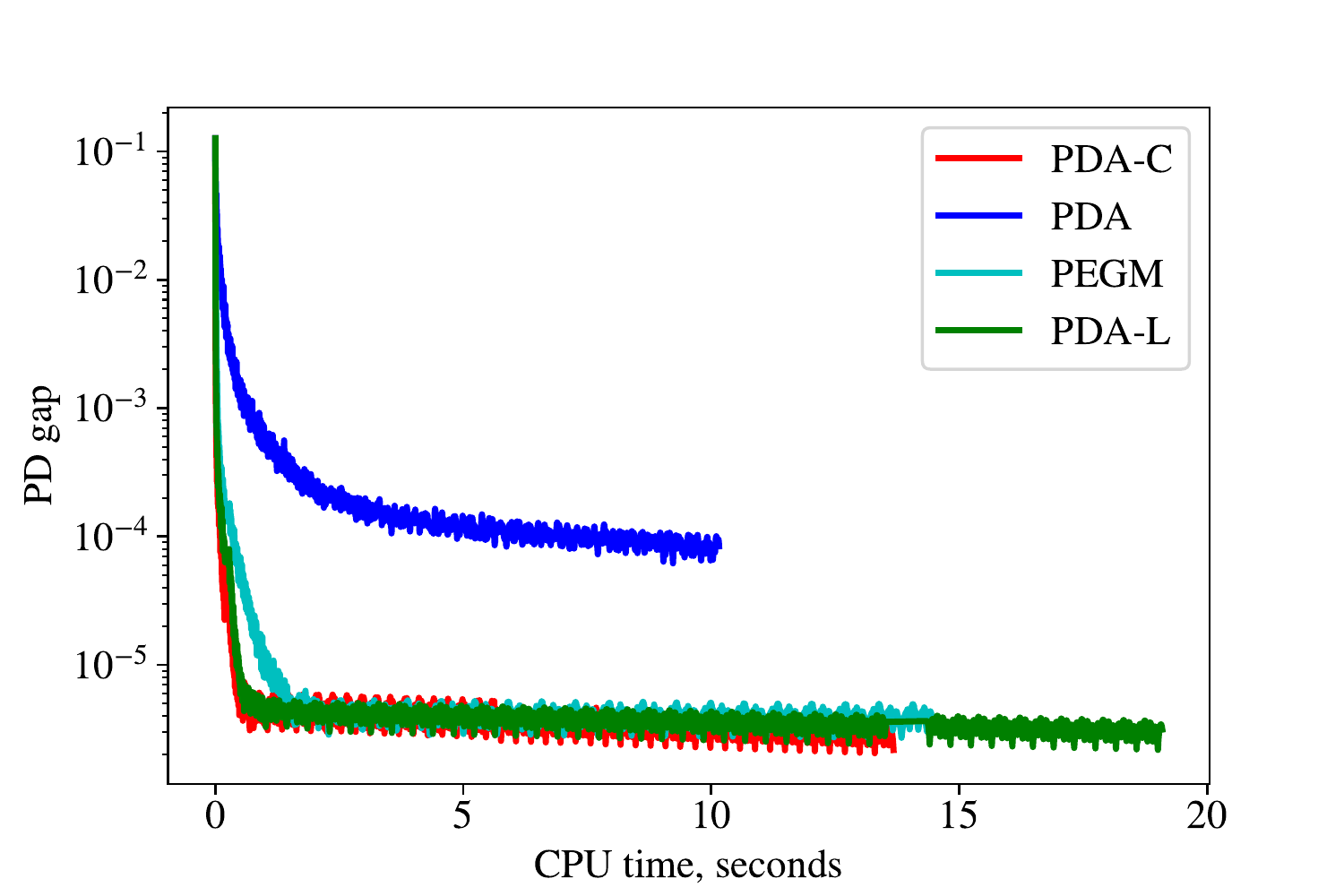}}
\caption{Comparison of PD gap for solving Problem \ref{pro_3} within 100,000 iterations. }
\label{Fig 5} 
\end{figure}

For every case we report the computing time (Time) measured in seconds, and show the primal-dual gap vs the computing time. The results are presented in Fig. \ref{Fig 5}. The execution time of all iterations for PDA is the lowest, PDA-C is slightly more than PDA, and PDA-L is about 2 times more expensive than PDA. From Fig. \ref{Fig 5}, PDA-L and PDA-C show better performance than PDA for the instances of Problem \ref{pro_3}. Furthermore, we notice that PDA-C can be better than PDA-L.

\begin{prob}[Nonnegative least square problem]\label{pro_2}
We are interested in the following problem
\ben
\min_{x>0} \phi(x):=\frac 1 2 ||Kx-b||^2
\een
or in a primal-dual form
\be\label{pro_22}
\min_{x\in \R^n} \max_{y\in \R^m} g(x)+\langle Kx,y\rangle-f^*(y),
\ee
where $f(p) = \frac 1 2 ||p-b||^2$, $f^*(y) = \frac 1 2 ||y||^2 + (b,y) = \frac 1 2 ||y+b||^2 -\frac{1}{2}||b||^2$, $g(x) = \delta_{\R^n_+}(x)$.
\end{prob}

We consider two real data examples from the Matrix Market library \footnote{https://math.nist.gov/MatrixMarket/data/Harwell-Boeing/lsq/lsq.html.} One is ``WELL1033": sparse matrix with $m = 1033, n = 320$, another is ``ILLC1033": sparse matrix with $m = 1033, n = 320$. For all cases the entries of vector $b\in \bR^m$ are generated independently from $\cN (0, 1)$.

To apply FISTA, we define $h(x) = f(Kx) = \frac{1}{2} ||Kx-b||^2$ for Problem \ref{pro_2}, then $\nabla h(x) = K^*(Kx-b)$. Since $f^*$ is strongly convex, in addition to PDA, PDA-L, and FISTA, we include in our comparison  APDA, APDA-L and APDA-C. We apply APDA-L to the primal-dual form (\ref{pro_22}) and APDA-C to the symmetry of (\ref{pro_22}), namely,
\ben
f^*(y) = \delta_{\R^n_+}(y),~~~g(x) =\frac 1 2 ||x+b||^2 -\frac{1}{2}||b||^2.
\een
We take parameter of strong convexity as $\gamma=1/2$. For PDA, APDA, and FISTA we compute $\|K\|$ and set $\tau= \sigma = \frac{1}{\|K\|}$, $\alpha=\frac{1}{\|K\|^2}$.
For PDA-L and PDA-C we set $\beta= 1$ (the same as $\tau=\sigma$ in PDA) and for PDA-C we set $\hat{n}=5000$. Since $\lambda_n\rightarrow0$ from Lemma \ref{acc_beta}, we set $\phi_n\equiv1$ for APDA-C in this section. The initial points are $x_0 = (0,\cdots, 0)$ and $y_0 = Kx_0-b =-b$.

\begin{figure}[htp]
\centering
\subfigure[$\phi(x)-\phi^*$]{
\includegraphics[width=0.45\textwidth]{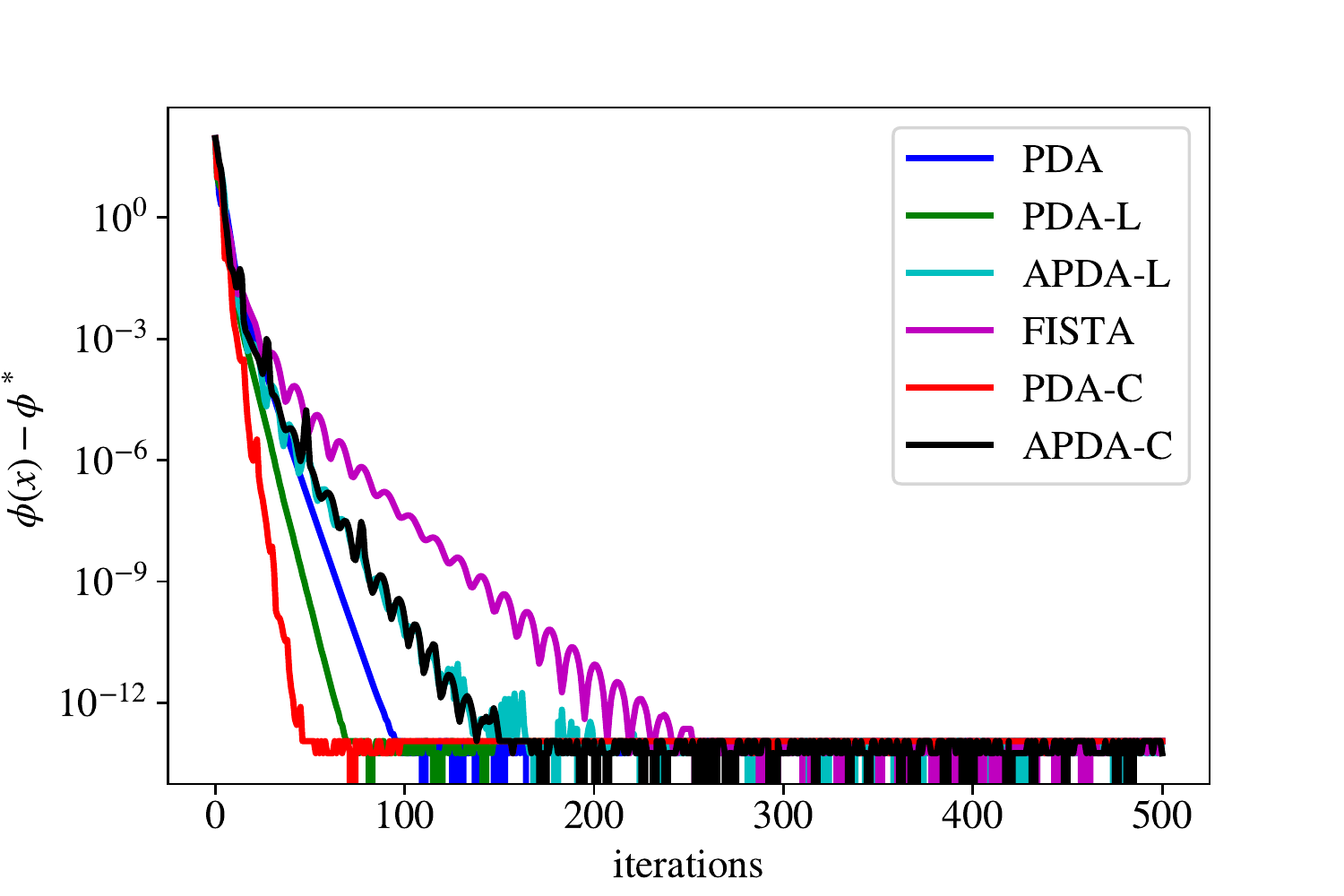}}
\subfigure[$\lambda_n$(or $\beta_{n-1}\lambda_n$, $\tau_k$)]{
\includegraphics[width=0.45\textwidth]{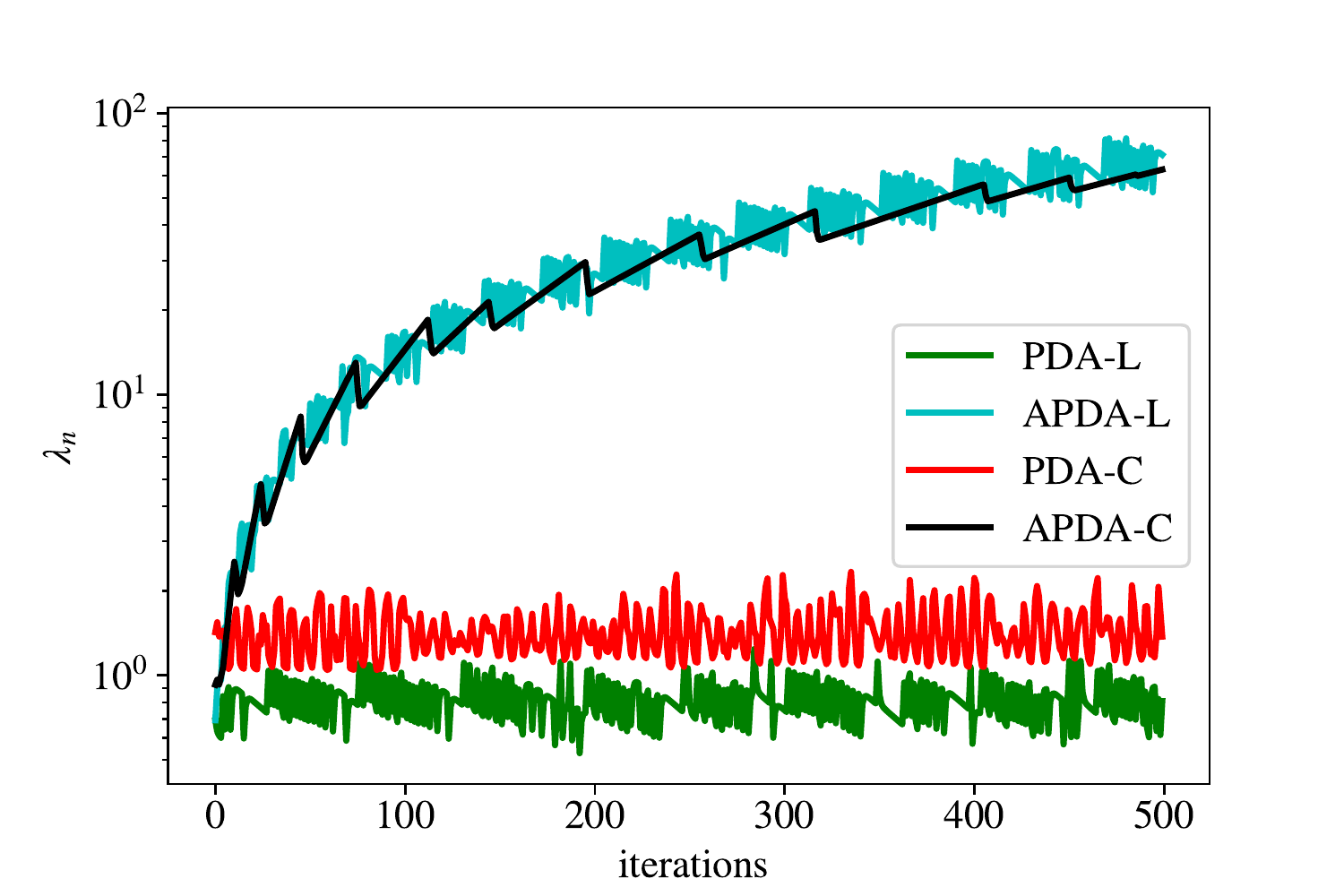}}
\caption{Comparison of $\phi(x)-\phi^*$ and $\lambda_n$(or $\beta_{n-1}\lambda_n$, $\tau_k$) for solving Problem \ref{pro_2} with ``WELL1033". }
\label{Fig 3}
\end{figure}

\begin{figure}[htp]
\centering
\subfigure[$\phi(x)-\phi^*$]{
\includegraphics[width=0.45\textwidth]{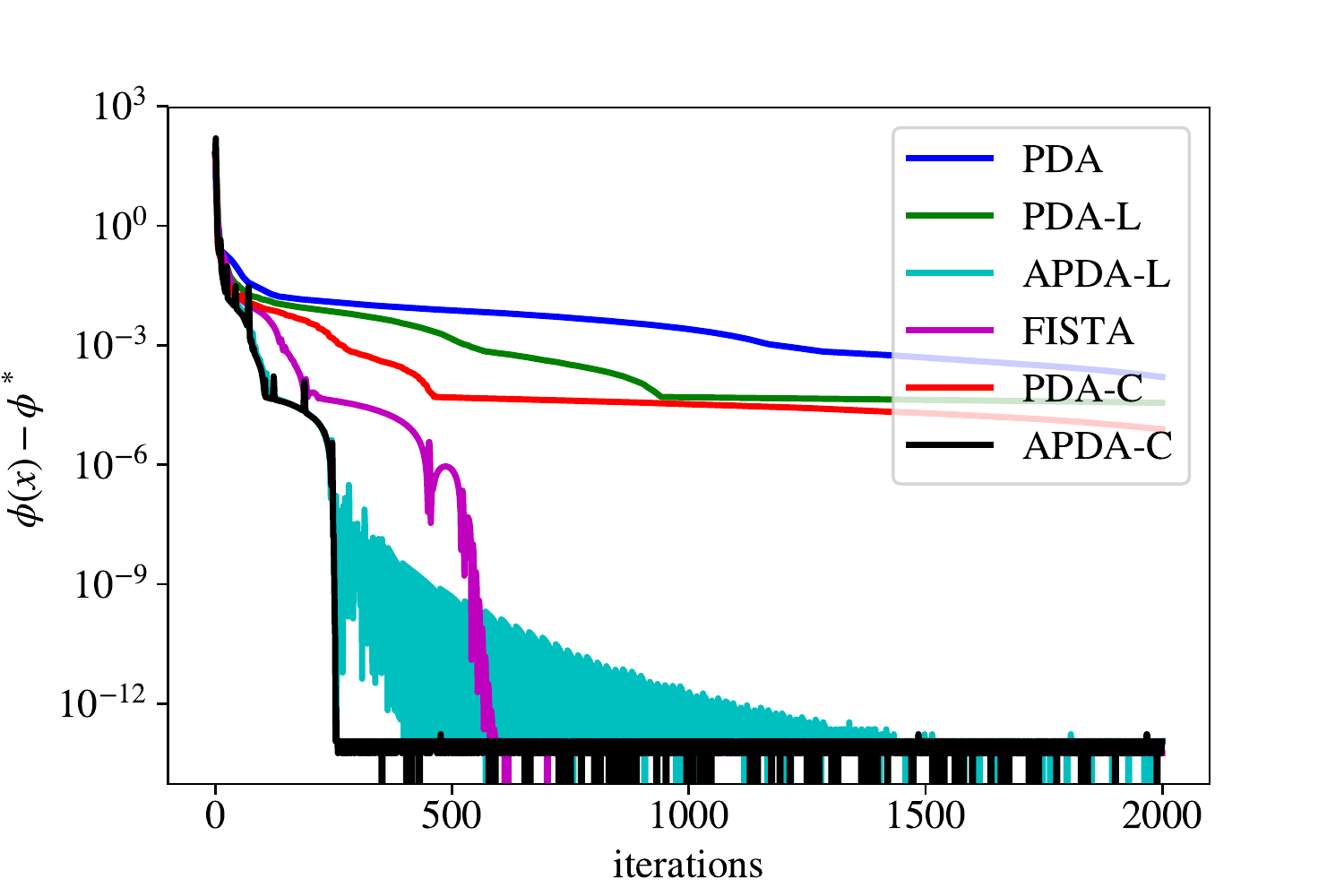}}
\subfigure[$\lambda_n$(or $\beta_{n-1}\lambda_n$, $\tau_k$).]{
\includegraphics[width=0.45\textwidth]{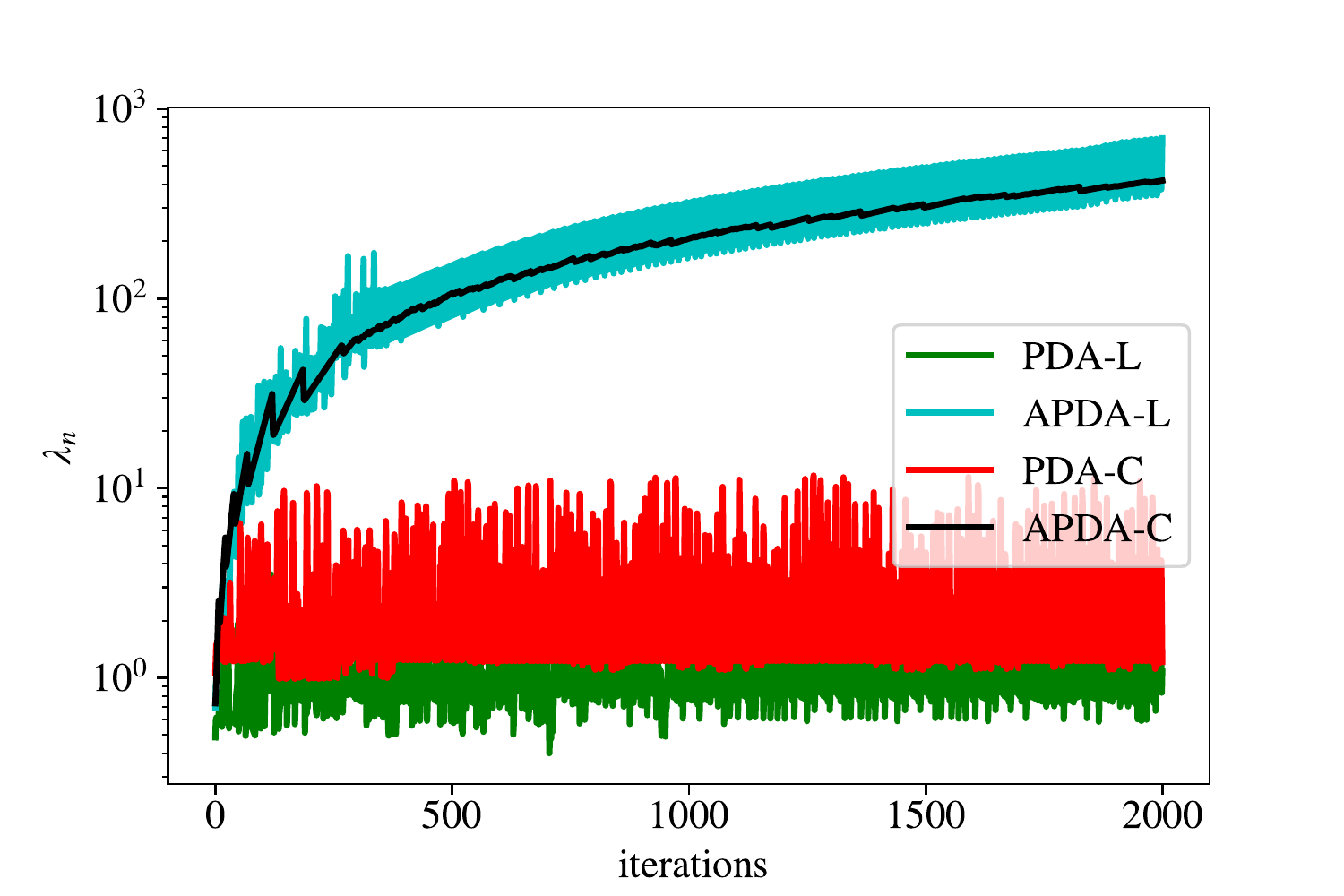}}
\caption{Comparison of $\phi(x)-\phi^*$ and $\lambda_n$(or $\beta_{n-1}\lambda_n$, $\tau_k$) for solving Problem \ref{pro_2} with ``ILL1033". }
\label{Fig 4} 
\end{figure}

We illustrate the plots of  $\phi(x_n)-\phi^*$ and $\lambda_n$ from PDA-C (or $\beta_{n-1}\lambda_n$ from APDA-C, $\tau_k$ from PDA-L) change over iterations. From the results shown in Fig. \ref{Fig 3} and \ref{Fig 4}, we observe that PDA-C with $\delta=0.62$ and PDA-L show better performance for the case ``WELL1033", while APDA-C and APDA-L for the case ``ILL1033". It is interesting to highlight that sometimes non-accelerated methods can be better than their accelerated variants for the well problems.

\section{Conclusions}
\label{sec_conclusion}
In this work, we have presented a primal-dual algorithm with correction and explored its acceleration. Firstly, the proposed PDA-C allows us to avoid the evaluation of the operator norm. Secondly, we have presented a prediction-correction strategy to estimate step sizes, which may results in larger step sizes. In practice, the correction step is conducted infrequently as only a weak condition needs to satisfy. Finally, we have proved convergence and established convergence rate, for PDA-C and its accelerated version.

Notice that only a very weak condition needs to check in the correction, and the correction step is not necessary to arrive termination conditions for some problems. Whether do the proposed PDA converge without correction? For $\delta\in]0,\frac{\sqrt{5}-1}{2}]$, whether there are any (larger) $\alpha>0$ such that the proposed PDA is convergent? We leave these as interesting topics for our future research.

\begin{acknowledgements}
The research of Xiaokai Chang was supported by  the Hongliu Foundation of First-class Disciplines of Lanzhou University of Technology. The project was supported by  the National Natural Science Foundation of China under Grant 61877046.
\end{acknowledgements}


\end{document}